\documentclass[reqno,11pt]{article}

\usepackage{natbib}
\bibpunct{(}{)}{,}{a}{}{;}

\usepackage{dsfont}

\usepackage[T1]{fontenc}      % Police contenant les caracteres francais
\usepackage{ae}                      % Pour ne pas que ce soit baveux
\usepackage[latin1]{inputenc}
\usepackage{lmodern}
\usepackage{amsmath}
\usepackage{amsfonts}
\usepackage[french,english]{babel}
\usepackage{amssymb}
\usepackage{dsfont}
\usepackage{array}
\usepackage{color}
\usepackage{epsfig}
\usepackage{fancybox}
\usepackage{amsthm}
\usepackage{float}
\usepackage{shorttoc}
\usepackage{color}
\usepackage[colorlinks=true,citecolor=blue]{hyperref}
\usepackage{shorttoc}
\usepackage[boxruled]{algorithm2e}
\usepackage{pdfsync}
\usepackage{enumitem}
\usepackage[colorlinks=true,citecolor=blue]{hyperref}
 \usepackage{slashbox}
\usepackage{caption}
\usepackage{subcaption}
\usepackage{amsthm, amsmath, amsfonts, amssymb}%
\usepackage{appendix} 
\usepackage{chngcntr} 
\usepackage{etoolbox} 
\usepackage{lipsum} 
\usepackage{rotating}
\usepackage[normalem]{ulem}
\usepackage{pifont}
\usepackage{float}
\usepackage{multirow}

\usepackage[T1]{fontenc}      % Police contenant les caracteres francais
\usepackage{ae}                      % Pour ne pas que ce soit baveux
\usepackage[latin1]{inputenc}
\usepackage{lmodern}
\usepackage{amsmath}
\usepackage{amsfonts}
\usepackage[french,english]{babel}
\usepackage{amssymb}
\usepackage{dsfont}
\usepackage{array}
\usepackage{color}
\usepackage{epsfig}
\usepackage{fancybox}
\usepackage{amsthm}
\usepackage{frcursive}
\usepackage{float}
\usepackage{shorttoc}

\theoremstyle{plain} 
\newtheorem{theorem}{Theorem}[section]
\newtheorem{lemma}[theorem]{Lemma} 
\newtheorem{corollary}[theorem]{Corollary} 
\newtheorem{proposition}[theorem]{Proposition} 
\newtheorem{assumption}[theorem]{Assumption}{\bf}{\rm}%
\newtheorem{remark}[theorem]{Remark}

\theoremstyle{remark}

\DeclareMathOperator{\pen}{pen}

\DeclareMathOperator{\Card}{Card}
\DeclareMathOperator{\Sp}{Sp}

\DeclareMathOperator{\var}{Var}

\DeclareMathOperator{\MISE}{MISE}

\DeclareMathOperator{\MISErand}{MISErand}

\DeclareMathOperator{\ISErand}{ISErand}

\newcommand{\bs}{\boldsymbol}
\newcommand{\sumin}{\displaystyle\sum_{i=1}^n}

\newcommand{\inttau}{\displaystyle\int_{0}^{\tau}}

\newcommand{\e}{{\mathrm e}}
\newcommand{\diff}{\mathrm d}

\newcommand{\hatbeta}{\bs{\hat\beta}}

\newcommand{\exphatbeta}{\e^{\bs{\hat\beta^TZ_i}}}
\newcommand{\expbeta}{\e^{\bs{\beta_0^TZ_i}}}

\newcommand{\rand}{rand(\bs{\hat\beta})}

\newcommand{\E}{\mathbb{E}}
\newcommand{\Proba}{\mathbb{P}}

\newcommand{\alphainf}{||\alpha_0||_{\infty,\tau}}

\usepackage{geometry}
\usepackage{pdfsync}

\RequirePackage{amsthm, amsmath, amsfonts, amssymb}%

\usepackage{graphicx, color}%
\usepackage{tikz}             % Pour des dessins
\usepackage{mathrsfs}

\usepackage[boxruled]{algorithm2e}

\definecolor{darkblue}{rgb}{0.0,0.0,0.7}

\hypersetup{%
  pdfauthor = {Agathe Guilloux, Sarah Lemler and Marie-Luce Taupin},%
  pdftitle = {},%
  pdfcreator = {pdflatex},%
  pdfproducer = {pdflatex}}

\date{}

\begin{document}

%\begin{frontmatter}
\title{Adaptive estimation of the baseline hazard function in the Cox model by model selection, with high-dimensional covariates}
%\begin{aug}
%\author{Agathe Guilloux
%\ead[label=e1]{agathe.guilloux@upmc.fr}
%\ead[label=u1,url]{http://www.lsta.lab.upmc.fr}}
%\address{Laboratoire de Statistique Th\'eorique et Appliqu\'ee\\
%Universit\'e Pierre et Marie Curie - Paris 6, \\
%4, place Jussieu\\
%75005 Paris, FRANCE\\
% \printead{u1}
%}
%
% \author{Sarah Lemler\ead[label=e2]{sarah.lemler@genopole.cnrs.fr}}
% \and
%  \author{Marie-Luce Taupin
%  \ead[label=e3]{marie-luce.taupin@genopole.cnrs.fr}\
%  \ead[label=u2,url]{http://stat.genopole.cnrs.fr}}
%
%  \address{Laboratoire Statistique et G\'enome\\
%UMR CNRS 8071, INRA 1152, Universit\'e d'\'Evry\\
%523, place des Terrasses de l'Agora\\
%91000 \'Evry, FRANCE\\
%          \printead{e2,e3,u2}}
%
%%  \thankstext{t2}{Footnote to the first author with the `thankstext' command.}
%
%  \runauthor{A. Guilloux et al.}
%
%\end{aug}
%  $^{(1)}${\footnote{(1) LSTA University Pierre et Marie Curie, France}}, Sarah Lemler $^{(1,2)}$, Marie-Luce Taupin $^{(2)}${\footnote{(2) Laboratoire Statistique et G\'enome, Universit\'e d'\'Evry Val d'Essonne, France}}}
%\address{(1) Laboratoire de Statistique Th\'eorique et Appliqu\'ee}
%\\
%University Pierre et Marie Curie\\
%%%...\\
%%%...\\
%Paris}

%\address{(2) Laboratoire Statistique et G\'enome\\
%Universit\'e d'\'Evry Val d'Essonne\\
%UMR CNRS 8071- USC INRA \\
%23 Boulevard de France\\
%91037 \'Evry}

%\maketitle
\author{Agathe Guilloux\\
\small{
Laboratoire de Statistique Th\'eorique et Appliqu\'ee,}\\
\small{Universit\'e Pierre et Marie Curie - Paris 6}\\
\small{\textit{e-mail :} \texttt{agathe.guilloux@upmc.fr}}\\
%\vspace{0.5cm}
\and
Sarah Lemler\\
\small{
Laboratoire de Math\'ematiques et de Mod\'elisation d'Evry, UMR CNRS 8071- USC INRA,}\\
\small{Universit\'e d'\'Evry Val d'Essonne, France}\\
\small{\textit{e-mail :} \texttt{sarah.lemler@genopole.cnrs.fr}}\\
\and
Marie-Luce Taupin\\
\small{Laboratoire de Math\'ematiques et de Mod\'elisation d'Evry, UMR CNRS 8071- USC INRA,}\\
\small{Universit\'e d'\'Evry Val d'Essonne, France}\\
\small{\textit{e-mail :} \texttt{marie-luce.taupin@genopole.cnrs.fr}}\\
\small{Unité MaIAGE, INRA Jouy-En-Josas, France}
}

%\author{Sarah Lemler\\
%\small{
%Laboratoire Statistique et G\'enome UMR CNRS 8071- USC INRA,}\\
%\small{Universit\'e d'Evry Val d'Essonne, France}\\
%\small{\textit{e-mail :} \texttt{sarah.lemler@genopole.cnrs.fr}}
%}

\maketitle

\begin{abstract}
The purpose of this article is to provide an adaptive estimator of the baseline function in the Cox model with high-dimensional covariates. We consider a two-step procedure : first, we estimate the regression parameter of the Cox model via a Lasso procedure based on the partial log-likelihood, secondly, we plug this Lasso estimator into a least-squares type criterion and then perform a model selection procedure to obtain an adaptive penalized contrast estimator of the baseline function. 

Using non-asymptotic estimation results stated for the Lasso estimator of the regression parameter, we establish a non-asymptotic oracle inequality for this penalized contrast estimator of the baseline function, which highlights the discrepancy of the rate of convergence when the dimension of the covariates increases.
\begin{flushleft}
\end{flushleft}
\end{abstract}

%\begin{keyword}[class=MSC]
%\kwd[Primary ]{60K35}
%\kwd{60K35}
%\kwd[; secondary ]{60K35}
%\end{keyword}
%
%\begin{keyword}
%\kwd{Survival analysis; Conditional hazard rate function; Cox's proportional hazards model; Right-censored data; Semi-parametric model; Nonparametric model; High-dimensional covariates; Model selection; Non-asymptotic oracle inequalities; concentration inequalities.}
%%\kwd{\LaTeXe}
%\end{keyword}

%\tableofcontents

%\end{frontmatter}

\section{Introduction}
\label{Intro}

Consider the following Cox model, introduced by \citet{Cox} and defined, for a vector of covariates $\bs{Z}=(Z_1,...,Z_p)^T$, by
\begin{align}
\label{eq:CoxIntro34}
\lambda_0(t, \bs{Z}) = \alpha_0(t)\exp(\bs{\beta_0^TZ}),
\end{align}
where  $\lambda_0$ denotes the hazard rate, $\bs{\beta_0} = (\beta_{0_1} , ..., \beta_{0_p})^T \in \mathbb{R}^p$ is the regression parameter and $\alpha_0$ is the baseline hazard function. 
%Both parameters  are unknown and have to be estimated. 
The Cox partial log-likelihood, introduced by \citet{Cox}, allows to estimate $\bs{\beta_0}$ without the knowledge of $\alpha_0$, considered as a functional nuisance parameter. For the estimation of $\alpha_0$, one common way is to use a two step procedure, starting with the estimation of $\bs{\beta_0}$ alone and then
to plug this estimator into a non parametric type estimator $\alpha_0$, usually a kernel type estimator.

%Although there are two unknown parameters $\alpha_0$ and $\bs{\beta_0}$  in the Cox model, more attention has been given to the estimation of the regression parameter $\bs{\beta_0}$, considering $\alpha_0$ as a functional nuisance parameter.
%One of the reasons is that  the Cox partial log-likelihood, introduced by \citet{Cox}, allows to estimate $\bs{\beta_0}$ without the knowledge of $\alpha_0$.
%%While it is true that we need to estimate both parameters in the Cox model to predict the survival time, in practice, the estimation of the baseline function has its own interest. 
%
%For the estimation of $\alpha_0$, one common way is to use a two step procedure, starting with the estimation of $\bs{\beta_0}$ alone and then
%to plug this estimator into a non parametric type estimator $\alpha_0$, usually a kernel type estimator. 
%%Those two steps would be changed according to the order of $p$ compared to the sample size $n$. 
%Very few results exists on the estimation of the baseline function $\alpha_0$. 
 
 Let us be more specific. 
 
When $p$ is small  compared to $n$, $\bs{\beta_0}$ is usually estimated by minimization of the opposite of the Cox partial log-likelihood.
 We refer to \citet{ABG}, as a reference book, for the proofs of the consistency and the asymptotic normality of $\bs{\hat\beta}$ when $p$ is small compared to $n$.
Thoses strategies only apply when $p<n$  and even more, they only apply when $p$ is small compared to $n$.
When $p$ growths up, becoming of the same order as $n$ and possibly larger than $n$, various well known problems appears. Among them, the minimization of the opposite of  the Cox partial log-likelihood
becomes difficult and even impossible if $p>n$.

  In high-dimension, when $p$ is large compared to $n$, the Lasso procedure is one of the classical  considered strategies. The Lasso (Least Absolute Shrinkage and Selection Operator)
  has been first introduced by \citet{Tib} in the linear regression model. It has been largely considered in additive regression model (see for instance \citet{knight00}, \citet{efron04}, \citet{donoho06},  \citet{Meinshausen}, \citet{Zhao06}, \citet{Zhang1}, \citet{Meinshausen09} and also \citet{juditsky00}, \citet{nemirovski00}, \citet{BTW2,BTW1,Bunea07}, \citet{greenshtein04} or
   \citet{BRT}), and  in density estimation (see  \citet{bunea2007sparse} and \citet{BlPR}). In the particular case of the semi-parametric Cox model,  \citet{T} has proposed a Lasso procedure for the regression parameter.
The Lasso estimator of the regression parameter  $\bs{\hat\beta}$ is defined as the minimizer of the opposite of the   Cox partial log-likelihood under an  $\ell_1$ type constraint, that is, suitably penalized with an $\ell_1$-penalty function. Recent results exist on the estimation of $\bs{\beta_0}$ 
in high-dimension setting. Among them one can mention \citet{Fan} who have proved asymptotic results for Lasso estimator. More recently, \citet{bradic2012}, \citet{Kong} and \citet{Huang2013} establish
the first non-asymptotic oracle inequalities (estimation and prediction bounds) for the Lasso estimator. 

For the baseline hazard function and when $p$ is small compared to $n$, the common estimator is a kernel estimator, which depends on $\bs{\hat\beta}$ obtained by minimization of the opposite of the Cox partial log-likelihood. This kernel estimator has been introduced by \citet{ramlau83b,ramlau83a} from the Breslow estimator of the cumulative baseline function (see \citet{ramlau83a} and \citet{ABG} for more details). In this context, \citet{ramlau83a} and \citet{Gregoire93} proved asymptotic results. No non-asymptotic results and no adaptive results have to date been established for the kernel estimator of the baseline function. Finally, when $p$ is large compared to $n$,  to our knowledge, the construction of an estimator of the baseline function has not been yet considered.

 In this paper, we consider a two-step procedure to estimate  $\bs{\beta_0}$ and $\alpha_0$, the two parameters  in the Cox model.
%  but the step isbased on model selection procedure. 
But our contributions focus more on  the estimation of $\alpha_0$.
%via a model selection procedure. 
In the Cox model we consider, it is noteworthy that the high-dimension only concerns the regression parameter, whereas the baseline function is a time function. 
Its estimation would not require a procedure specific to high-dimension, besides the  first step concerning the estimation of $\bs{\beta_0}$.
We propose a procedure for the construction of  an estimator of the baseline hazard function $\alpha_0$,  $p$ being either smaller than $n$ or greater than $n$. It combines a Lasso procedure for  $\bs{\beta_0}$ as a first step and a second step based on a model selection strategy for the estimation of the baseline function $\alpha_0$. 
%\textcolor{red}{As mentionned above, the Lasso estimator of $\bs{\beta_0}$ is obtained by minimization of  the opposite of the Cox partial log-likelihood 
%penalized by an $\ell_1$-penalty whereas the estimator of $\alpha_0$ is obtained by penalized least squares criterion minimisation over models collection. }
This model selection procedure takes its origins in the works of \citet{Akaike} and \citet{Mallows73}, more recently formalized by \citet{birge97} and \citet{barron99} for the estimation of densities and regression functions (see the book of \citet{Massart} as a reference work on model selection). In survival analysis, the model selection has also been documented. \citet{Letue} has adapted these methods to estimate the regression function of the non-parametric Cox model, when $p<n$.
More recently, \citet{brunel05}, \citet{brunel09}, \citet{BCL10} have obtained adaptive estimation of densities in a censoring setting. Model selection methods have also been used to estimate the intensity function of a counting process in the multiplicative Aalen intensity model (see \citet{Reynaud06} and \citet{CGG}). However, the model selection procedure has never been considered, to our knowledge, for estimating the baseline hazard function in the Cox model.

Our contributions are at least threefold:
 Our procedure  is the first that focus on the estimation of baseline function of the semi-parametric Cox model with high-dimentional covariates.
This procedure provide an adaptive estimator of the baseline function that works as well for small $p$ and large $p$ compared to $n$ (that is for possibly high-dimensional covariates). 
 Furthermore, for this estimator, we  state non-asymptotic oracle inequalities, that hold, once again, $p$ being either smaller than $n$ or greater than $n$. More precisely, we prove that the risk of this estimator achieves the best risk among estimators in a large collection. For each model, the risk of an estimator is bounded by the sum of three terms. The first term is a bias term involving to the approximation properties of the collection of models, through the distance evaluated in  $\bs{\beta_0}$ between the true baseline and the orthogonal projection of $\alpha_0$ on the best selected model. 
The second term is a penalty term of the same order than the variance on one model, that is of order the dimension of one model over $n$, as expected with $\ell_0$-penalty. These two terms are the "usual" terms appearing in nonparametric estimation. It is noteworthy that these two terms  do not involve any quantity related to the risk of the Lasso estimator of  $\bs{\beta_0}$. The last term precisely comes from the properties of the Lasso estimator of  $\bs{\beta_0}$. This last term is of order $\log(n p)/n$, as expected for a Lasso estimator.

When $p$ is small, the third last term is of order $\log(n)/n$ and, the rate is governed by the first two terms. In that case,  the penalty term being of the same order than the variance over one model, we conclude that  the model selection procedure achieves the "expected rate" of order $n^{-2\gamma/(2\gamma+1)}$ when the baseline function belongs to a Besov space with smoothness parameter $\gamma$. This continues to hold when
 $p$ is of the same order than the sample size $n$. When $p$ is larger than $n$, that is in the so-called ultra-high dimension (see
\citet{Verzelen12}), the rate for estimating $\alpha_0$ is changed, and more precisely degraded  as a price to pay for being with high dimension covariates. 
This degradation follows accordingly to the order of $p$ compared to $n$.

The main tools for stating our results are the theory of marked counting processes and martingales with jumps, the theory of penalized minimum contrast estimators and
concentrations inequalities such as Talagrand inequality (see \citet{Talagrand1996}) and a Bernstein inequality found in   (see \citet{SVG} and \citet{CGG})  for unbounded martingale process and combined with chaining methods (see \citet{Talagrand05generic} and \citet{Baraud10bernstein}).

The article is organized as follows. In Section \ref{sec:description2}, we describe the estimation procedure. Section \ref{sec:NAOI2} provides non-asymptotic oracle inequalities on the estimator of the baseline hazard function $\alpha_0$, in a high-dimensional setting for $\bs{\beta_0}$. In section \ref{sec:simu}, we compare the performances of the resulting penalized contrast estimator to those of the usual kernel estimator on simulated data. Section \ref{sec:technic2} is devoted to the proofs: we 
state some technical results, then we establish the two main theorems and lastly we prove the technical results. 
%The performances of the resulting penalized contrast estimator to those of the usual kernel estimator on simulated data. Finally, Section \ref{sec:technic2} is devoted to the proofs: we 
Finally, Appendix \ref{ann:Huang} discusses the bound of the error estimation for the Lasso estimator of the regression parameter of the Cox model.

%%% SECTION 2 --- DESCRIPTION OF THE PROCEDURE ---
\section{Notations and preliminaries}

\subsection{Framework with counting processes}
\label{sub:framework}
Consider the general setting of counting processes, which embeds the classical case of right censoring. We follow here the now classical setting of \citet{ABG} or \citet{fleming11}. For $n$ independant individuals, we observe for $i=1,...,n$ a counting process $N_{i}$, a random process $Y_{i}$ with values in $[0,1]$ and a vector of covariates $\bs{Z_i}=(Z_{i,1},...,Z_{i,p})^T\in\mathbb{R}^p$. Let $(\Omega,\mathcal{F},\mathbb{P})$ be a probability space and $(\mathcal{F}_{t})_{ t\geq 0}$ be the filtration 
defined by 
$$\mathcal{F}_{t}=\sigma\{N_{i}(s),Y_{i}(s), 0\leq s\leq t, \bs{Z_{i}}, i=1,...,n\}.$$
From the Doob-Meyer decomposition, we know that each $N_i$ admits a compensator denote by $\Lambda_i$, such that $M_{i}=N_{i}-\Lambda_{i}$ is a $(\mathcal{F}_{t})_{t\geq 0}$  local square-integrable martingale (see \citet{ABG} for details). We assume in the following that $N_i$ has a satisfies an Aalen multiplicative intensity model. 
\begin{assumption}
For each $i=1,...,n$ and all $t\geq 0$, 
\begin{align}
\label{eq:aalenIntro34}
\Lambda_i(t)=\displaystyle{\int_{0}^{t}\lambda_{0}(s,\bs{Z_{i}})Y_{i}(s)\diff s},
\end{align}
where $\lambda_0(t,\bs{z})=\alpha_0(t)\e^{\bs{\beta^Tz}}$, for $\bs{z}\in\mathbb{R}^p$.
\end{assumption}

We observe the independent and identically distributed (i.i.d.) data $(\bs{Z_{i}}, N_{i}(t), Y_{i}(t), i=1,...,n, 0\leq t\leq \tau)$, where $[0,\tau]$ is the time interval between the beginning and the end of the study. \\

This general setting, introduced by \citet{aalen}, embeds several particular examples as censored data, marked Poisson processes and Markov processes (see \citet{ABG} for further details). 
%The intensity $\lambda_0$ has to be estimated.
We give here details for the right censoring case. We observe for $i=1,...,n$, $(X_i,\delta_i,\bs{Z_i})$, where $X_i=\min({T_{i}},{C_{i}})$, $\delta_{i}=\mathds{1}_{\{T_{i}\leq C_{i}\}}$, $T_i$ is the time of interest and $C_i$ the censoring time. With these notations, the $(\mathcal{F}_{t})_{t\geq 0}$-adapted processes $Y_{i}$ and $N_{i}$ are respectively defined as the at-risk process $Y_{i}(t)=\mathds{1}_{\{X_{i}\geq t\}}$ and the counting process $N_{i} (t)=\mathds{1}_{\{X_{i}\leq t,\delta_{i}=1\}}$ which jumps when the ith individual dies.

\subsection{Assumptions}

Before describing the estimation procedure, we introduce few assumptions on the framework defined in Subsection \ref{sub:framework}. 

Let $\bs{Z}\in\mathbb{R}^p$ denote the generic vector of covariates with the same distribution as the vectors of covariates $\bs{Z_i}$ of each individual $i$ and by $Z_j$ its   $j$-th component, namely the $j$-th covariates of the vector $\bs{Z}$. Similarly, we denote by $Y$ the generic version of the random process $Y_i$ with values in $[0,1]$.

We define the standard $\mathbb{L}^{2}$ and $\mathbb{L}^{\infty}$-norms, for $\alpha\in (\mathbb{L}^{2}\cap \mathbb{L}^{\infty})([0,\tau])$:
\[
||\alpha||^2_2=\inttau \alpha^2(t)\diff t \quad \mbox{and} \quad ||\alpha||_{\infty,\tau}=\underset{t\in[0,\tau]}{\sup}|\alpha(t)|.
\]
For a vector $\bs{b}\in\mathbb{R}^p$, we also introduce the $\ell_1$-norm $|b|_1=\sum_{j=1}^{p} |b_j|$.
\newpage
%Let us introduce some assumptions on the covariates $(\bs{Z}_i)_{1\leq i\leq n}$ and on $\alpha_0$.

\begin{assumption}\
\label{ass:assnot2}
\begin{enumerate}[label=\textbf{(\roman*)},ref=(\roman*)]

\item \label{ass:Zij2} There exists a positive constant $B$ such that  
\[
|Z_{j}|\leq B, \quad\forall j\in\{1,...,p\}.
\]
In the following, we denote $A=[-B,B]^p$.

\item \label{ass:f1} The vector of covariates $\bs{Z}$ admit a p.d.f. $f_{\bs Z}$ such that $\sup_ A|f_{\bs Z}|\leq f_1<+\infty$.

\item \label{ass:f0} There exists $f_0>0$, such that $\forall (t,\bs{z})\in [0,\tau]\times A$, 
\[
\mathbb{E}[Y(t)|\bs{Z}=\bs{z}]f_{\bs{Z}}(\bs{z})\geq f_0.
\]

\item \label{ass:alpha0inf2} For all $t\in[0,\tau]$, $\alpha_0(t)\leq ||\alpha_0||_{\infty,\tau}<+\infty$.\\

\end{enumerate}
\end{assumption}

\begin{remark} \
Let say a few word on these assumptions starting by noting that these four assumptions are quite classic and reasonnable. To be more specific, Assumption \ref{ass:assnot2}.\ref{ass:Zij2}, is  very common to establish oracle inequalities of Lasso estimators in various frameworks.
In particular,   in the Cox model,  see e.g. \citet{Huang2013} and \citet{bradic2012} for the statement of non asymptotic oracle inequalities

In the specific case of right censoring, Assumption \ref{ass:assnot2}.\ref{ass:f0} is automatically verified. Indeed, for $T$ the survival time and $C$ the censoring time, we can write 
\[\mathbb{E}(Y(t)|\bs{Z}=\bs{z})=\mathbb{E}(\mathds{1}_{\{T\wedge C\leq t\}}|\bs{Z}=\bs{z})=(1-F_{T|\bs{Z}}(t))(1-G_{C|\bs{Z}}(t-)),
\]
where $F_{T|\bs{Z}}$ and $G_{C|\bs{Z}}$ are the cumulative distribution functions of $T|\bs{Z}$ and $C|\bs{Z}$ respectively.
It is known (see \citet{ABG}) that the Kaplan-Meier estimator is consistent only on intervals of the form $[0,\tau]$, 
where $\tau\leq \sup\{t\geq 0, (1-F_{T|\bs{Z}}(t))(1-G_{C|\bs{Z}}(t))>0\}$. Hence when $f_{\bs{Z}}$ is bounded from below on $A$, there exists $f_0>0$, such that 
\begin{equation*}
\forall (t,\bs{z})\in[0,\tau]\times A, \quad \mathbb{E}[Y(t)|\bs{Z}=\bs{z}]f_{\bs{Z}}(\bs{z})\geq f_0.
\end{equation*}

Assumption \ref{ass:assnot2}.\ref{ass:f0}  is required in order to compare the natural norm of the baseline function induced by our contrast to the standard $\mathbb{L}^2$-norm (see Proposition \ref{prop:connection2}).

\end{remark}

\section{Estimation procedure}
\label{sec:description2}
%We are interested in obtaining an adaptive estimator of the baseline function via model selection, in a high-dimensional setting for $\bs{\beta_0}$. This estimation relies on a two-step procedure. 
We now describe our two-steps estimation procedure, starting by recalling the Lasso estimation of $\bs{\beta_0}$ and then giving a bound of its prediction risk. Then, we describe the contrast and the model selection procedure for the estimation of the baseline function.

\subsection{Preliminary estimation of $\bs{\beta_0}$: procedure and results}
\label{sub:estimationbeta}

%As usual and recalled in the introduction of this chapter, the baseline function relies on the preliminary estimation of the regression parameter $\bs{\beta_0}$. 
%Let us describe the first preliminary step of our procedure on the estimation of the regression parameter $\bs{\beta_0}$ in high-dimension. 
%We consider a classical Lasso procedure defined in Section \ref{Intro} by (\ref{estimationbeta0Intro34}) and (\ref{eq:penIntro34}) to estimate $\bs{\beta_0}$.
The Lasso estimator $\bs{\hat\beta}$ of the regression parameter $\bs{\beta_0}$, introduced in \citet{T}, is defined by
% is estimated via a classical Lasso procedure defined by: 
%The Cox partial log-likelihood has been introduced by \citet{Cox} to estimate the regression parameter $\bs{\beta_0}$ in the Cox model without the knowledge of the baseline function. It is 
%defined, for all $\bs\beta\in\mathbb{R}^p$, by 
%\begin{equation}
%\label{eq:PartialCoxLikelihood2}
% l^*_n(\bs\beta)=\dfrac{1}{n}\sumin\inttau\log\dfrac{\e^{\bs{\beta^TZi}}}{S_n(t,\bs{\beta})}\diff N_i(t), \quad \mbox{where } S_n(t,\bs\beta)=\dfrac{1}{n}\sumin \e^{\bs{\beta^TZ_i}}Y_i(t).
%\end{equation}
%From this partial log-likelihood, we define the usual estimator of $\bs{\beta_0}$ when $p<n$ and the Lasso estimator of $\bs{\beta_0}$ when $p\gg n$: 
\begin{align}
\label{eq:betaL2}
\bs{\hat\beta}=\underset{\bs\beta\in\mathbb{R}^p}{\arg\min}\{- l^*_n(\bs{\beta})+\Gamma_n|\bs{\beta}|_1\}, 
%\label{eq:Lasso2}
% \mbox{with} \quad \pen(\bs{\beta})&=\left\{
%   \begin{array}{ll}
%       0 & \mbox{when } p<n \\
%        \mu|\bs{\beta}|_1 & \mbox{in high-dimension when } p\gg n, 
%   \end{array}
%\right.
\end{align}
 where $\Gamma_n$ is a positive regularization parameter to be suitable chosen, $|\bs\beta|_1=\sum_{j=1}^p|\beta_j|$ and $ l^*_n$ is the Cox partial log-likelihood defined by, 
\begin{equation}
\label{eq:PartialCoxLikelihood2}
 l^*_n(\bs\beta)=\dfrac{1}{n}\sumin\inttau\log\dfrac{\e^{\bs{\beta^TZi}}}{S_n(t,\bs{\beta})}\diff N_i(t), \quad \mbox{where } S_n(t,\bs\beta)=\dfrac{1}{n}\sumin \e^{\bs{\beta^TZ_i}}Y_i(t) \quad \forall t\geq 0.
\end{equation}

The risk bounds for the estimator of $\alpha_0$ will naturally involve the risk $|\bs{\hat\beta}-\bs{\beta_0}|_1$,  that have to be at least bounded. Thus, we rather consider the following procedure 
\begin{equation}
\label{estb0Ball2}
\bs{\hat\beta}=\underset{\bs\beta\in\mathcal{B}(0,R_1)}{\arg\min}\{- l^*_n(\bs{\beta})+\pen(\bs\beta)\}, \quad \mbox{with} \quad \pen(\bs{\beta})= \Gamma_n|\bs{\beta}|_1,
\end{equation}
where $\mathcal{B}(0,R_1)$ is the ball defined by
\[
\mathcal{B}(0,R_1)=\{b\in\mathbb{R}^p : |b|_1\leq R_1\}, \quad \mbox{with } R_1>0.
\]
Consider the following assumption:

%We now describe assumptions on $\bs{\hat\beta}$ and $\bs{\beta_0}$,
%%This estimator $\bs{\hat\beta}$ and the true parameter $\bs{\beta_0}$ satisfy some assumptions, that are 
%required to ensure the existence of the estimator of the baseline function defined in the second step.
\begin{assumption}\
\label{ass:betaball2}
%\begin{enumerate}[label=\textbf{(\roman*)},ref=(\roman*)]
%\item \label{ass:hatbeta2} $\bs{\hat\beta}\in\mathcal{B}(0,R_1)$, where $\mathcal{B}(0,R_1)$ is the ball defined by
%\[
%\mathcal{B}(0,R_1)=\{b\in\mathbb{R}^p : |b|_1\leq R_1\}, \quad \mbox{with } R_1>0.
%\]
%On this ball, the procedure becomes
%\begin{equation}
%\label{estb0Ball2}
%\bs{\hat\beta}=\underset{\bs\beta\in\mathcal{B}(0,R_1)}{\arg\min}\{- l^*_n(\bs{\beta})+\pen(\bs\beta)\}, \quad \mbox{with} \quad \pen(\bs{\beta})= \Gamma_n|\bs{\beta}|_1,
%\end{equation}
%\item \label{ass:beta02} $|\bs{\beta_0}|_1<R_2<+\infty$.
We assume that $|\bs{\beta_0}|_1<R_2<+\infty$.

\end{assumption}
We denote $R=\max(R_1,R_2)$, so that 
\begin{equation}
\label{eq:beta0-betaR2}
|\bs{\hat\beta}-\bs{\beta_0}|_1\leq 2R \quad \mbox{ a.s}.
\end{equation}
Such condition has already been considered by \citet{SVG1} or \citet{Kong}. Roughly speaking, it means that we can restrict our attention to a ball, possibly very large, in a neighborhood of $\bs{\beta_0}$ for finding a good estimator of $\bs{\beta_0}$. 

As mentionned above, our risk bounds for the estimator of $\alpha_0$ depend on the risk $|\bs{\hat\beta}-\bs{\beta_0}|_1$.  Such bounds on this risk already exist. In particular, in their Theorem 3.1, \citet{Huang2013} state a non asymptotic inequality for $|\bs{\hat\beta}-\bs{\beta_0}|_1$ in the specific case of bounded counting processes. We consider here more general processes, possibly unbounded. 
In the following proposition, we provide a generalization of the results established by \citet{Huang2013} to the case of unbounded counting processes.
We refer to Appendix \ref{ann:Huang} for a proof of Proposition \ref{prop:predbeta2}.

\begin{proposition}
\label{prop:predbeta2}
Let $k>0$, $c>0$ and $s:=\Card\{j\in\{1,...,p\}: \beta_{0_j}\neq 0\}$ be the sparsity index of $\bs{\beta_0}$. Assume that $||\alpha_0||_{\infty,\tau}<\infty$. Then, under Assumptions \ref{ass:betaball2} and \ref{ass:Zij2}, with probability larger than $1-cn^{-k}$, we have
\begin{equation}
\label{predictionbeta}
|\bs{\hat\beta}-\bs{\beta_0}|_1\leq C(s)\sqrt{\dfrac{\log (pn^k)}{n}}
\end{equation}
where $C(s)>0$ is a constant depending on the sparsity index $s$.
\end{proposition}
%\begin{remark}
%Under Assumption \ref{ass:baseline3}.\ref{ass:alpha0inf3}, Lemma \ref{lem:OmegaH} in Chapter \ref{chap:ModelSelection} Subsection \ref{sub:estimationbeta} is directly proved from Lemma 3.3 from \citet{Huang2013}.
%\end{remark}
As mentioned previously, this proposition is crucial to establish a non-asymptotic oracle inequality for the baseline function. 
In the rest of the paper, we consider that  $\bs{\hat\beta}$ satisfies Inequality (\ref{predictionbeta}).

\begin{assumption}
\label{ass:Relationpn}
We assume that 
\[
\underset{n\rightarrow\infty}\lim C(s)\dfrac{\log(np)}{n}=0.
\]
\end{assumption}
This assumption is clearly reasonable: when $p$ is smaller than $n$ or of the same order, this assumption is automatically fulfilled.
It is not satisfied 
when $p$ becomes too high compared to $n$. This case corresponds to the now well known case of ultra-high dimension framework. In this specific case, recent lower bounds in additive regression models typically say that the estimation of paramater is mostly impossible (see for example \citet{Verzelen12}).

\subsection{Estimation of $\alpha_0$}
We now come to the estimation of the baseline function $\alpha_0$ via a model selection procedure. As usual, such a procedure requires an empirical estimation criterion, a collection of models and a suitable penalty function, all being presented in the following.

\subsubsection{Definition of the estimation criterion}
\addcontentsline{minitoc}{subsubsection}{Definition of the estimation criterion}
We estimate the baseline function $\alpha_0$ using a least-squares criterion. More precisely, based on the data $(\bs{Z_{i}}, N_{i}(t), Y_{i}(t), i=~1,...,n, 0\leq ~t\leq \tau)$ and for a fixed $\bs\beta$, 
we consider the empirical least-squares type given for a function $\alpha\in (\mathbb{L}^{2}\cap \mathbb{L}^{\infty})([0,\tau])$ by  
\begin{equation}
\label{eq:contrast2}
C_{n}(\alpha,\bs{\beta})=-\dfrac{2}{n}\displaystyle{\sum_{i=1}^{n}\int_{0}^{\tau}\alpha(t)\diff N_i(t)+\dfrac{1}{n}\sum_{i=1}^{n}\inttau\alpha^2(t)\e^{\bs{\beta^TZ_i}}Y_i(t)\diff t}.
\end{equation}
%It is not very common to use least-squares criterion in survival analysis. 
The use of such 
least-square empirical criterion in survival analysis is not so usual as for the additive regression model. Nevertheless, few recent studies have developped such very useful as strategies. Among them one can cite  \citet{Reynaud06} or \citet{CGG}.

Let us define a deterministic scalar product and its associated deterministic norm for $\alpha_1$, $\alpha_2$ and $\alpha$ functions in $(\mathbb{L}^{2}\cap \mathbb{L}^{\infty})([0,\tau])$:
\begin{eqnarray}
\nonumber
\langle\alpha_1,\alpha_2\rangle_{det(\bs\beta)}\hspace{0.2cm}&=&\inttau\alpha_1(t)\alpha_2(t)\mathbb{E}[e^{\bs{\beta^TZ}}Y(t)] \diff t,\\
\label{normdet}
||\alpha||_{det(\bs\beta)}^2&=&
\inttau\alpha^2(t)\mathbb{E}[e^{\bs{\beta^TZ}}Y(t)]\diff t. \end{eqnarray}

Using the Doob-Meyer decomposition $N_i=M_i+\Lambda_i$ and according to the multiplicative Aalen model (\ref{eq:aalenIntro34}), we get:
\[
\mathbb{E}[C_n(\alpha,\bs{\beta_0})]=||\alpha||^2_{det}-2\langle \alpha,\alpha_0 \rangle_{det}=||\alpha-\alpha_0||^2_{det}-||\alpha_0||^2_{det},
\]
which is minimum when $\alpha=\alpha_0$.
Hence, 
%the estimation criterion (\ref{eq:contrast2}) is suitable and 
minimizing $C_n(., \bs{\beta_0})$ is a relevant strategy to estimate $\alpha_0$.

\subsubsection{Model selection}
\label{subsec:estproc}
We now describe the model selection procedure in our context, introducing first the collection of models.

\paragraph{Collections of models.}
Let $\mathcal{M}_n$ be a set of indices and $\{S_m,m\in\mathcal{M}_n\}$ be a collection of models:
\[
S_m=\{\alpha : \alpha=\sum_{j\in J_m}a_j^m\varphi_j^m, a_j^m\in\mathbb{R}\},
\]
where $(\varphi_j^m)_{j\in J_m}$ is an orthonormal basis of $(\mathbb{L}^2\cap \mathbb{L}^{\infty})([0,\tau])$ for the usual $\mathbb{L}_2(P)$- norm.
We denote $D_m$ the cardinality of $S_m$, i.e. $|J_m|=D_m$.

\paragraph{Sequence of estimators.}
Let us consider $\bs{\hat\beta}$ the Lasso estimator of $\bs{\beta_0}$ defined by (\ref{estb0Ball2}).
For each $m\in\mathcal{M}_n$, we define the estimator
 \begin{align}
 \label{eq:hatalphamChap3}
 \hat\alpha_m^{\bs{\hat\beta}}=\underset{\alpha\in S_m}{\arg\min }\{ C_n(\alpha,\bs{\hat\beta})\}.
\end{align}

\paragraph{Model selection.}
The relevant space is automatically selected by using following penalized criterion
\begin{align}
\label{eq:mhathatbetaChap3}
\hat m^{\bs{\hat\beta}}=\underset{m\in\mathcal{M}_n}{\arg\min}\{C_n(\hat\alpha^{\bs{\hat\beta}}_m,\bs{\hat\beta})+\pen(m)\},
\end{align}
where $\pen: \mathcal{M}_n\rightarrow\mathbb{R}$ will be defined later.

\paragraph{Final estimator.}
The final estimator of $\alpha_0$ is then $\hat\alpha^{\bs{\hat\beta}}_{\hat m^{\bs{\hat\beta}}}$.
\newline

Let us say few words on the optimisation problem.
Denote by  $\bs{G^{\bs{\hat\beta}}_m}$  the random Gram matrix
\begin{equation}
\label{eq:gram2}
\bs{G^{\bs{\hat\beta}}_m}=\Big(\dfrac{1}{n}\sumin\inttau\varphi_j(t)\varphi_k(t)\e^{\bs{\hat\beta^TZ_i}} Y_i(t)\diff t\Big)_{(j,k)\in J^2_m}.
\end{equation}
By definition, the estimator $\hat\alpha_m^{\bs{\hat\beta}}$ is the  solution of the equation $\bs{G^{\bs{\hat\beta}}_mA^{\hat\beta}_m}=\bs{\Gamma_m}$, where  
\begin{equation}
\label{Gamma}
\bs{A^{\bs{\hat\beta}}_m}=(\hat{a}^{\bs{\hat\beta}}_j)_{j\in J_m} \quad \mbox{and} \quad\bs{\Gamma_m}=\Big(\dfrac{1}{n}\sumin\inttau\varphi_j(t)\diff N_i(t)\Big)_{j\in J_m}.
\end{equation}
The Gram matrix $\bs{G^{\bs{\hat\beta}}_m}$ may not be invertible in some cases.  Hence we consider the set 
\begin{align}
\label{eq:Hm2}
\hat{\mathcal{H}}^{\bs{\hat\beta}}_m=\Bigg\{\min\Sp(\bs{G^{\bs{\hat\beta}}_m})\geq \max\Bigg(\dfrac{\hat{f}_0\e^{-B|\bs\beta_0|_1}\e^{-B|\bs\beta_0-\bs{\hat\beta|_1}}}{6},\dfrac{1}{\sqrt{n}}\Bigg)\Bigg\},
\end{align}
where $\Sp(\bs{M})$ denotes the spectrum of matrix $\bs{M}$ and $\hat f_0$ satisfies the following assumption:

\begin{assumption}
\label{ass:hatf0} 
There exist a preliminary estimator $\hat f_0$ of $f_0$ and two positive constants $C_0>0$, $n_0>0$ such that 
%Let denote $\hat f_0$ the estimator of $f_0$. 
%there are positive constants $C_0$ and $n_0$ such that 
\[
\mathbb{P}(|\hat f_0-f_0|>f_0/2)\leq C_0/n^6 \quad \mbox{for any} \quad n\geq n_0.
\] 
\end{assumption}

 From Assumptions \ref{ass:betaball2}, on the set $\hat{\mathcal{H}}^{\bs{\hat\beta}}_m $, the matrix $\bs{G^{\bs{\hat\beta}}_m}$ is invertible and $\hat\alpha^{\bs{\hat\beta}}_m$ is thus uniquely defined as 
\[
\hat\alpha^{\bs{\hat\beta}}_m=\left\{
    \begin{array}{ll}
    	\arg\min_{\alpha\in\mathcal{S}_m}\{ C_n(\alpha,\bs{\hat\beta})\} & \mbox{on } \hat{\mathcal{H}}^{\bs{\hat\beta}}_m,\\
	0 & \mbox{on } (\hat{\mathcal{H}}^{\bs{\hat\beta}}_m)^c.
    \end{array}
 \right.		
\]

\subsubsection{Assumptions and examples of the models}

The following assumptions on the models $\{S_m : m\in\mathcal{M}_n\}$ are usual in model selection procedures. They are verified by the spaces spanned by usual bases: trigonometric basis, regular  piecewise polynomial basis, regular compactly supported wavelet basis and histogram basis. We refer to \citet{barron99} and \citet{brunel05} for other examples and further discussions.
\begin{assumption}\
\label{ass:model2}
\begin{enumerate}[label=\textbf{(\roman*)},ref=(\roman*)]
\item \label{ass:Dm2}
For all $m\in\mathcal{M}_n$, we assume that
\[
D_m\leq \dfrac{\sqrt{n}}{\log n}.
\]
\item \label{ass:norminf2}
For all $m\in\mathcal{M}_n$, there exists $\phi>0$ such that for all $\alpha$ in $S_m$,
\[
\underset{t\in[0,\tau]}{\sup}|\alpha(t)|^2\leq \phi D_m\inttau\alpha^2(t)\diff t.
\]

\item \label{ass:emboite} The models are nested within each other: $D_{m_1}\leq D_{m_2} \Rightarrow S_{m_1}\subset S_{m_2}$.
We denote by $\mathcal{S}_n$ the global nesting space in the collection and by $\mathcal{D}_n$ its dimension. 

\end{enumerate}
\end{assumption}
\begin{remark}
Assumption \ref{ass:model2}.\ref{ass:Dm2} ensures that the sizes $D_m$ of the models are not too large compared with the number of observations $n$. This assumption seems reasonable if we remember that $D_m$ is the number of coefficients to be estimated: if this number is too large compared to the size of the panel, we cannot expect to obtain a relevant estimator.
Assumption \ref{ass:model2}.\ref{ass:norminf2} implies a useful connection between the standard $\mathbb{L}^{2}$-norm and the infinite norm. Assumption \ref{ass:model2}.\ref{ass:emboite} ensures that $\forall m,m'\in\mathcal{M}_n$, $S_m+~S_{m'}\subset \mathcal{S}_n$. Thanks to this assumption, one does not have to browse through all models for the model selection, which reduces the algorithmic complexity of the procedure.
In addition, we have from Assumption \ref{ass:model2}.\ref{ass:Dm2} that $\mathcal{D}_n\leq \sqrt{n}/{\log n}$.
\end{remark}

%%% SECTION 3 --- NAOI ON THE BASELINE ---

\section{Non-asymptotic oracle inequalities}
\label{sec:NAOI2}

We now are in a position to state our main theorem: a non-asymptotic oracle inequality for the estimator $\hat\alpha^{\bs{\hat\beta}}_{\hat m^{\bs{\hat\beta}}}$ of the baseline function in the Cox model.
%\textcolor{blue}{
%Let $\alpha^{\bs{\beta_0}}_m$ be the projection of $\alpha_0$ on $S_m$ with respect to the deterministic scalar product when $\bs{\beta_0}$ is known:
%\begin{align}
%\label{def:alpham}
%\alpha^{\bs{\beta_0}}_m=\underset{\alpha\in S_m}{\arg\min }\hspace{0.1cm}\mathbb{E} [C_n(\alpha,\bs\beta_0)]=\underset{\alpha\in S_m}{\arg\min }||\alpha-\alpha_0||^2_{det}.
%\end{align}
%The estimator $\hat\alpha^{\bs{\hat\beta}}_{\hat m^{\bs{\hat\beta}}}$ satisfies the following inequality:
%}

\begin{theorem}
\label{th:IO2}
Let Assumptions \ref{ass:assnot2}.\ref{ass:Zij2}-\ref{ass:alpha0inf2}, Assumptions \ref{ass:betaball2}, Assumption \ref{ass:Relationpn}, Assumption \ref{ass:hatf0} and Assumptions \ref{ass:model2}.\ref{ass:Dm2}-\ref{ass:emboite} hold. Let $\alpha^{\bs{\beta_0}}_m$ be the projection of $\alpha_0$ on $S_m$ with respect to the deterministic scalar product when $\bs{\beta_0}$ is known:
\begin{align}
\label{def:alpham}
\alpha^{\bs{\beta_0}}_m=\underset{\alpha\in S_m}{\arg\min }\hspace{0.1cm}\mathbb{E} [C_n(\alpha,\bs\beta_0)]=\underset{\alpha\in S_m}{\arg\min }||\alpha-\alpha_0||^2_{det}.
\end{align}
Let $\hat\alpha^{\bs{\hat\beta}}_{\hat m^{\bs{\hat\beta}}}$  be defined by (\ref{eq:hatalphamChap3}) and (\ref{eq:mhathatbetaChap3}) with
\begin{equation}
\label{eq:pen2}
\pen(m):=K_0(1+||\alpha_0||_{\infty,\tau})\dfrac{D_m}{n},
\end{equation}
where $K_0$ is a numerical constant. Then, for any $n\geq n_0$, with $n_0$ a constant defined in Assumption \ref{ass:hatf0},
\begin{equation}
\label{eq:NAOI2}
\mathbb{E}[||\hat\alpha^{\bs{\hat\beta}}_{\hat m^{\bs{\hat\beta}}}-\alpha_0||^2_{det}]\leq \kappa_0\underset{m\in\mathcal{M}_n}{\inf}\{||\alpha_0-\alpha^{\bs{\beta_0}}_m||^2_{det}+2\pen(m)\}+\dfrac{C_1}{n}+C_2C(s)\dfrac{\log(np)}{n},
\end{equation}
where $\kappa_0$ is a numerical constant, $C_1$ and $C_2$ are constants depending on $\tau$, $\phi$, $||\alpha_0||_{\infty,\tau}$, $f_0$, $\mathbb{E}[\e^{\bs{\beta_0^TZ}}]$, $\mathbb{E}[\e^{2\bs{\beta_0^TZ}}]$, $\mathbb{E}[\e^{4\bs{\beta_0^TZ}}]$, $B$, $|\bs{\beta_0}|_1$, the sparsity index $s$ of $\bs{\beta_0}$ and $\kappa_b$ a constant from the B\"urkholder Inequality (see Theorem \ref{th:Burkholder}) and $C(s)$ the constant depending on the sparsity index of $\bs{\beta_0}$ in Proposition \ref{prop:predbeta2}.
\end{theorem}

Inequality (\ref{eq:NAOI2}) provides the first non-asymptotic oracle inequality for an estimator of the baseline function. This inequality warrants the performances of our estimator $\hat\alpha^{\bs{\hat\beta}}_{\hat m^{\bs{\hat\beta}}}$. 
We refer to Subsection \ref{proof:IO2} for precisions about $C_1$ and $C_2$. In Inequality (\ref{eq:NAOI2}), the risk is bounded by the sum of four terms.

The third term of order $1/n$ is negligible compared to the others. 
The first two terms are respectively the bias and the variance terms. 
The bias term, $||\alpha_0-\alpha^{\bs{\beta}_0}_m||^2_{det}$, corresponds to the approximation error and decreases with the dimension $D_m$ of the model $S_m$. It depends on the regularity of the true function, which is unknown: the more regular $\alpha_0$ is, the smaller the bias is. The variance term $\pen(m)$ quantifies the estimation error and in contrary to the bias term, increases with $D_m$. It is of order $D_m/n$, which corresponds to the order of the  variance term on one model. These three first terms do not involve quantities related to the estimation error of the Lasso estimator of $\bs{\beta_0}$.

The last term precisely comes from the non-asymptotic control of $|\bs{\hat\beta}-\bs{\beta_0}|_1$ given by Proposition \ref{prop:predbeta2}. Indeed, we can rewrite Inequality (\ref{eq:NAOI2}) before using the bound of control (\ref{predictionbeta}):
\begin{equation*}
\mathbb{E}[||\hat\alpha^{\bs{\hat\beta}}_{\hat m^{\bs{\hat\beta}}}-\alpha_0||^2_{det}]\leq \kappa_0\underset{m\in\mathcal{M}_n}{\inf}\{||\alpha_0-\alpha^{\bs{\beta_0}}_m||^2_{det}+2\pen(m)\}+\dfrac{C_1}{n}+C_2\E[|\bs{\hat\beta}-\bs{\beta_0}|_1^2].
\end{equation*}
This inequality makes clearer the role of the first step of the procedure in the control of the estimator $\hat\alpha^{\bs{\hat\beta}}_{\hat m^{\bs{\hat\beta}}}$ of the baseline function. The bound obtained for this control is of order $\log(np)/n$, which explains the order of the fourth term. This term quantifies the influence of the high dimension on the estimation of the baseline hazard function. For small $p$, we obtain the expected rate of convergence in the case of a purely non-parametric estimation, but when is larger than $n$, the rate of convergence of the inequality is degraded. This is the price to pay for dealing with covariates in high dimension.

\begin{corollary}
\label{cor:Minimax}
Assume that $\alpha_0$ belongs to the Besov space $\mathcal{B}^{\gamma}_{2,\infty}([0,\tau])$, with smoothness $\gamma$. Then, under the assumptions of Theorem \ref{th:IO2}, 
\[
\mathbb{E}[||\hat\alpha^{\bs{\hat\beta}}_{\hat m^{\bs{\hat\beta}}}-\alpha_0||^2_{2}]\leq\tilde Cn^{-\frac{2\gamma}{2\gamma+1}}+C_2C(s) \dfrac{\log(np)}{n},
\]
where $\tilde C$ and $C_2$ are constants depending on $\tau$, $\phi$, $||\alpha_0||_{\infty,\tau}$, $f_0$, $\mathbb{E}[\e^{\bs{\beta_0^TZ}}]$, $\mathbb{E}[\e^{2\bs{\beta_0^TZ}}]$,  $B$, $|\bs{\beta_0}|_1$, the sparsity index $s$ of $\bs{\beta_0}$ and $\kappa_b$ a constant from the B\"urkholder Inequality (see Theorem \ref{th:Burkholder}) and $C(s)$ the constant depending on the sparsity index of $\bs{\beta_0}$ from Proposition \ref{prop:predbeta2}.
\end{corollary}

%\begin{remark}
%\label{remark:minimax}
%Assume that $\alpha_0$ belongs to the anisotropic Besov space $\mathcal{B}^{\gamma}_{2,\infty}([0,\tau])$, with smoothness $\gamma$. 
%From the proof of Corollary 1 in \citet{CGG}, we deduce that 
%\begin{align*}
%\mathbb{E}[||\hat\alpha^{\bs{\hat\beta}}_{\hat m^{\bs{\hat\beta}}}-\alpha_0||^2_{2}]&\leq \dfrac{\e^{B|\bs{\beta_0}|_1}}{f_0}\mathbb{E}[||\hat\alpha^{\bs{\hat\beta}}_{\hat m^{\bs{\hat\beta}}}-\alpha_0||^2_{det}]\\
%& \leq\tilde C_1\underset{m\in\mathcal{M}_n}{\inf}\Big\{D_m^{-2\gamma}+\pen(m)\Big\}+\tilde C_2\dfrac{\log(np)}{n}
%\end{align*}
From \citet{Reynaud06}, we know that, for an intensity function without covariates in a Besov space with smoothness parameter $\gamma$, the minimax rate is $n^{-{2\gamma}/{(2\gamma+1)}}$. We infer that this would also be the optimal rate in our case when the term $\log(np)/n$ is negligible, namely when $p<n$. 
%When the term $\log(np)/n$ is negligible ($p\leq n$), we expect to obtain a rate of convergence of order $n^{-2\gamma/(2\gamma+1)}$, which is certainly minimax from the minimax bounds established by \citet{Reynaud06} for an intensity function without covariates and by \citet{CGG} for the intensity function with one covariates.
However, when the high-dimension $p\gg n$ is reached, the remaining term $\log(np)/n$ is not negligible anymore and there is a loss in the rate of convergence, which comes from the difficulty to estimate $\bs{\beta_0}$. 
%\end{remark}

%This non-asymptotic oracle inequality is the first theoretical adaptive result on the baseline function. Without the last term in $\log(np)/n$, Inequality (\ref{eq:NAOI2}) would correspond to usual oracle inequalities obtained in classical frameworks using model selection methods (see \citet{Massart} for densities, \citet{BCL10} for cumulative distribution functions or \citet{CGG} for intensity functions). However, since we first estimate the regression parameter $\bs{\beta_0}$ via a Lasso procedure and we then plug this estimator to estimate the baseline function via model selection, we recover a supplementary term of order $\log(np)/n$. The control of $|\bs{\hat\beta}-\bs{\beta_0}|_1$ given by Proposition \ref{prop:predbeta2} is of order $\log(np)/n$. This explains the order of our remaining term. The term in $\log(p)/n$ comes from the high-dimensional setting. This term is typical of high-dimension (see \citet{Lemler12}). Thus, although our oracle inequality converges more slowly than the oracle inequalities in classical frameworks, due to this term, it highlights nevertheless the newness of our approach: the estimation of the baseline function of the Cox model in a high-dimensional setting.
%

%\end{enumerate}
%%% SECTION 4 --- TECHNICAL RESULTS ---

\section{Applications: simulation study}
\label{sec:simu}

%\subsection{Simulation study}

The aim of this section is to illustrate the behavior of the penalized contrast estimator $\hat\alpha^{\bs{\hat\beta}}_{\hat m^{\bs{\hat\beta}}}$ of the baseline function in the case of right censoring and to compare it with the usual kernel estimator with a bandwidth selected by cross-validation introduced by \citet{ramlau83a}.

\subsection{Simulated data}

Let consider the Cox model (\ref{eq:CoxIntro34}) in the case of right censoring. We consider a cohort of size $n$ and $p$ covariates. In the simulation study, several choices of $n$ and $p$ have been considered. The sample size $n$ takes the values $n=200$ and $n=500$ and $p$ varies between $p=\sqrt{n}$, being $15$ and $22$ respectively and $p=n$, referred to as the high-dimension case. 
%\begin{itemize}
%\item[-] $p=\sqrt{n}$, being $15$ and $22$ respectively and referred to as the intermediate case,
%\item[-] $p=n$, referred to as the high-dimension case.
%\end{itemize}  

\begin{figure}
\centering
\includegraphics[scale=0.5]{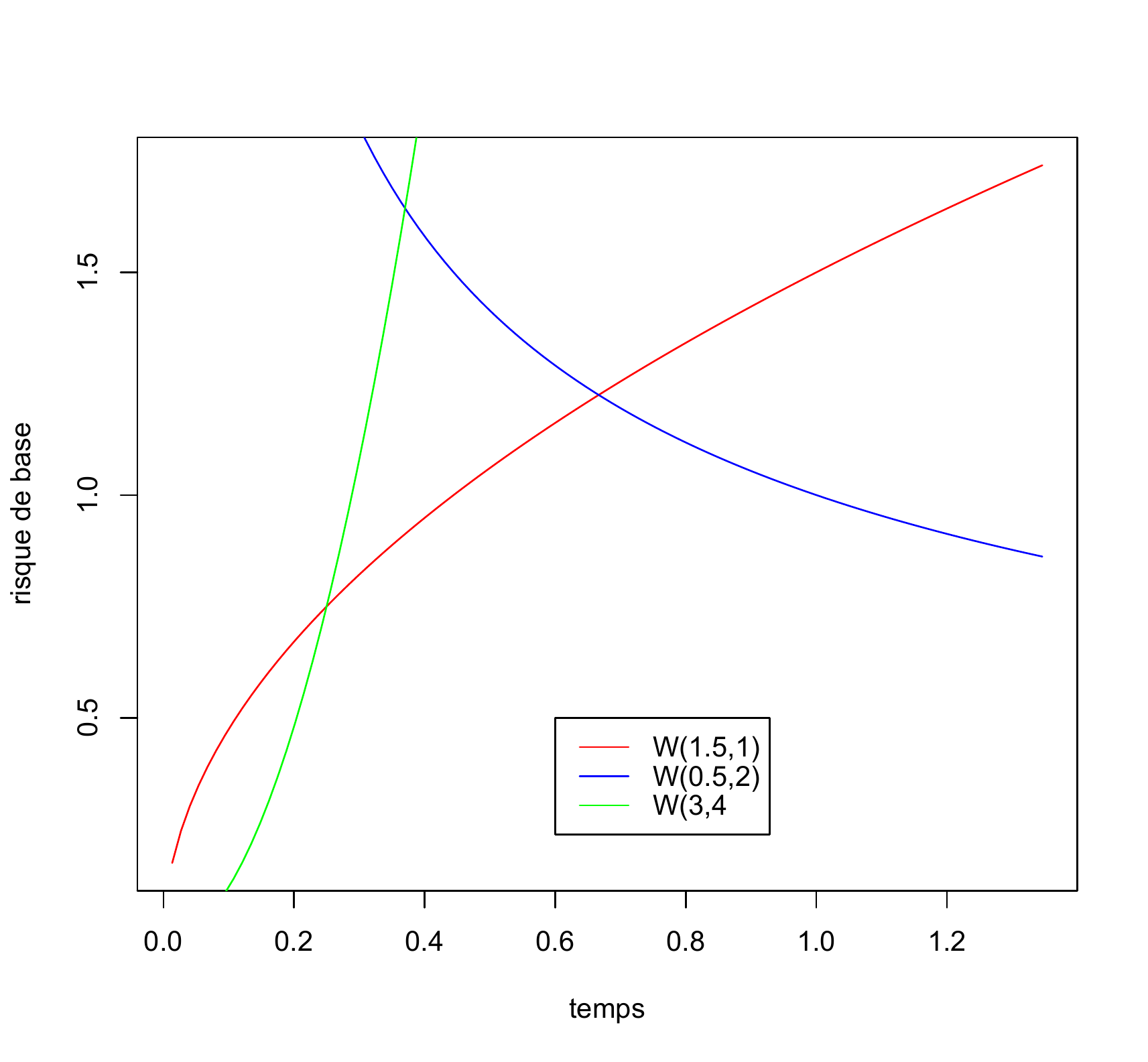}
\caption{Plots of the baseline hazard function for different parameters of a Weibull distribution $\mathcal{W}(a,\lambda)$}
\label{fig:baseline}
\end{figure}

The true regression parameter $\bs{\beta_0}$ is chosen as a vector of dimension $p$, defined by 
\[
\bs{\beta_0}=(0.1,0.3,0.5,0,...,0)^T\in\mathbb{R}^p,
\]
for various $p\geq3$ and for each $n$ and $p$, the design matrix $\bs{Z}=(Z_{i,j})_{1\leq i\leq n, 1\leq j\leq p}$ is simulated independently from a uniform distribution on $[-1,1]$. We consider survival times $T_i$, $i=1,...,n$ that are distributed according to a Weibull distribution $\mathcal{W}(a,\lambda)$, namely the associated baseline function is of the form $\alpha_0(t)=a\lambda^at^{a-1}$. 
We simulate three Weibull distribution $\mathcal{W}(0.5,1)$, $\mathcal{W}(1,1)$, $\mathcal{W}(3,4)$ (see Figure \ref{fig:baseline}). 
We consider a rate of censoring of $20\%$ and the censoring times $C_i$, for $i=1,...,n$,  are simulated independently from the survival times via an exponential distribution $\mathcal{E}(1/\gamma\E[T_1])$, where $\gamma=4.5$ is adjusted to the rate of censorship. The time $\tau$ of the end of the study is taken as the quantile at $90\%$ of $(T_i\wedge C_i)_{i=1,...,n}$. For $i=1,...,n$, we compute the observed times $X_i=\min(T_i,\tilde C_i)$, where $\tilde C_i=C_i\wedge\tau$ and the censoring indicators $\delta_i=\mathds{1}_{T_i\leq C_i}$. The definition of $\tilde C_i$ ensures that there exist some $i\in\{1,...,n\}$ for which $X_i\geq \tau$, so that all estimators are defined on the interval $[0,\tau]$ and it prevents from certain edge effect. 

Each sample $(\bs{Z_i}, T_i, C_i, X_i, \delta_i, i=1,...,n)$ is repeated $N_e=100$ times.

\subsection{Estimation procedures}

We implement $\hat\alpha^{\bs{\hat\beta}}_{ m}$ in a histogram basis defined, for $j=1,...,2^m$, by
\[
\varphi^m_j(t)=\dfrac{1}{\sqrt\tau}2^{m/2}\mathds{1}_{[{(j-1)\tau/2^m},{j\tau/2^m}[}(t),
\]
In this case, the cardinal of $S_m$ is $D_m=2^m$ and Assumption \ref{ass:model2}.\ref{ass:norminf2} is satisfied for $\phi=1/\tau$. We take $m=0,...,\lfloor\log(n/\log(n))/\log(2))\rfloor$, so that Assumption \ref{ass:model2}.\ref{ass:Dm2} is fulfilled.
In this basis, the estimator is being written by
\begin{equation}
\label{eq:BaseHist}
\hat\alpha^{\bs{\hat\beta}}_{m}(t)=\sum_{j\in J_{m}}\hat{a}^{\bs{\hat\beta}}_j\varphi^m_j(t), \quad \forall t\in[0,\tau],
\end{equation}
where 
\[
\hat a^{\bs{\hat\beta}}_j=\dfrac{\tau}{2^m}\dfrac{1}{\frac{1}{n}\sumin\e^{\bs{\hat\beta^TZ_i}}\Big(\Big(\min\Big(X_i,\frac{j\tau}{2^m}\Big)-\frac{(j-1)\tau}{2^m}\Big)\vee 0\Big)}\dfrac{1}{n}\sumin\delta_i\dfrac{2^{m/2}}{\sqrt{\tau}}\mathds{1}_{\Big[\frac{(j-1)\tau}{2^m},\frac{j\tau}{2^m}\Big)}(X_i).
\]
The final estimator $\hat\alpha^{\bs{\hat\beta}}_{\hat m^{\bs{\hat\beta}}}$ is obtained from the implementation of the selection model procedure  (\ref{eq:hatalphamChap3}), replacing in the penalty term the unknown quantity $\alphainf$ by $||\hat\alpha^{\bs{\hat\beta}}_{\max(m)}||_{\infty,\tau}$, an estimator of $\alpha_0$ computed on the arbitrary larger space $S_{\max(m)}$.

%We estimate the true regression parameter $\bs{\beta_0}$ by the Lasso procedure (\ref{eq:betaL2}). 
We want to compare the performances of the estimator $\hat\alpha^{\bs{\hat\beta}}_{\hat m^{\bs{\hat\beta}}}$ to those of the usual kernel estimator with a bandwidth selected by cross-validation introduced by \citet{ramlau83a}, that we have also implemented. More precisely the usual kernel estimator is defined by 
\begin{equation}
\label{eq:hatalphasimu4}
\hat\alpha^{\bs{\hat\beta}}_{\hat{h}^{\bs{\hat\beta}}_{CV}}(t)=\dfrac{1}{\hat{h}^{\bs{\hat\beta}}_{CV}}\sumin\dfrac{\delta_i}{\sum_{j=1}^{n}\e^{\bs{\hat\beta^TZ_j}}\mathds{1}_{\{X_j\geq X_i\}}}K\left(\dfrac{t-X_i}{\hat{h}^{\bs{\hat\beta}}_{CV}}\right),
\end{equation}
where $K(u)=0.75(1-u^2)\mathds{1}_{\{|u|\leq 1\}}$ is the Epanechnikov kernel and the bandwidth $\hat{h}^{\bs{\hat\beta}}_{CV}$ has been selected by cross-validation:
\[
\hat h^{\bs{\hat\beta}}_{CV}=\underset{h}{\arg\min}\Bigg\{\E\inttau(\hat\alpha^{\bs{\hat\beta}}_h(t))^2\diff t-2\sum_{i\neq j}\dfrac{1}{h}K\Big(\dfrac{X_i-X_j}{h}\Big)\dfrac{\Delta N(X_i)}{\bar Y(X_i)}\dfrac{\Delta N(X_j)}{\bar Y{(X_j)}}\Bigg\},
\] 
where $\bar Y=\sum_{i=1}^n\mathds{1}_{\{X_i\geq t\}}$.

Both estimators of the baseline hazard function are defined from the Lasso estimator $\bs{\hat\beta}$ of the regression parameter defined by (\ref{eq:betaL2}).

The performances of these two estimators are evaluated via a random Mean Integrated Squared Error ($\MISErand$) adapted to the Cox model and defined by $\MISErand(\alpha,\bs{\hat\beta})=\E[\ISErand(\alpha,\bs{\hat\beta})]$, where the expectation is taken on $(T_i, C_i, \bs{Z_i})$ and 
%there different Integrated Squared Errors defined for some function $\alpha\in\mathbb{L}^2([0,\tau])$ by  
%the standard $\ISE$, the random $\ISE$ more adapted to the Cox model and the total $\ISE$ are respectively defined by
\begin{align}
%\nonumber
%\ISEstand(\alpha)&=\inttau(\alpha(t)-\alpha_0(t))^2\diff t,\\
\label{eq:MISErand}
\ISErand(\alpha,\bs{\hat\beta})=\dfrac{1}{n}\sumin\int_0^{X_i}(\alpha(t)-\alpha_0(t))^2\e^{\bs{\hat\beta^TZ_i}}\diff t,
%\nonumber
%\ISEtotal(\alpha,\bs\beta)&=\dfrac{1}{n}\sumin\inttau(\alpha(t)\e^{\bs{\beta^TZ_i}}-\alpha_0(t)\e^{\bs{\beta_0^TZ_i}})^2\diff t.
%\ISEtotal(\alpha,\bs{\beta})&=\dfrac{1}{n}\sumin\inttau(\alpha(t)\e^{\bs{\beta^TZ_i}}-\alpha_0(t)\e^{\bs{\beta_0^TZ_i}})^2\diff t.
\end{align}
%The associated random Mean Integrated Squared Error is then defined by $\MISErand(\alpha)=\E[\ISErand(\alpha)]$, where the expectation is taken on $(T_i, C_i, \bs{Z_i})$ (for sake of simplicity, we write $\MISErand(\alpha)$ even if the $\MISErand$ depends on $\bs\beta$). 
We obtain an estimation of the $\MISErand$ by taking the empirical mean for $N_e=100$ replications.
\newline

In Table \ref{tab:MISEempLois}, we give the random $\MISE$ of the penalized contrast estimator and of the kernel estimator with a bandwidth selected by cross-validation for different distributions of the survival times.\\

\begin{table}[htp!]
\centering
%\begin{subtable}[b]{0.8\linewidth}
\begin{tabular}{|c|c|c|c||c|c||c|c|}

\hline
\multicolumn{2}{|c|}{\backslashbox{Dimensions}{Distributions}} &  \multicolumn{2}{|c||}{ $\mathcal{W}(1.5,1)$}  & \multicolumn{2}{|c||}{$\mathcal{W}(0.5,2)$ }  & \multicolumn{2}{|c|}{  $\mathcal{W}(3,4)$ } \\
\hline
  \multirow{2}{*}{$n=200$} & $p=15$   & 0.072  &  0.021  & 0.626 & 1.09 &  5.26&   8.48   \\
  \cline{2-8}
                                          & $p=200$  &    0.071 & 0.020 &   0.613     & 1.09      &    5.30   &     8.33        \\
 \hline
\multirow{2}{*}{$n=500$} & $p=22$     &0.055  &  0.009   & 0.401   &   1.06 & 5.24&  7.48 \\
\cline{2-8}
                                        & $p=500$    &    0.059   &  0.008 & 0.402  &       1.06       &5.25   &   8.10         \\              
                                           \hline        
%    \multicolumn{2}{c|}{} & $\hat\alpha^{\bs{\hat\beta}}_{\hat m^{\bs{\hat\beta}}}$ &   $\hat\alpha^{\bs{\hat\beta}}_{\hat{h}^{\bs{\hat\beta}}_{CV}}$ &$\hat\alpha^{\bs{\hat\beta}}_{\hat m^{\bs{\hat\beta}}}$ &  $ \hat\alpha^{\bs{\hat\beta}}_{\hat{h}^{\bs{\hat\beta}}_{CV}}$ & $\hat\alpha^{\bs{\hat\beta}}_{\hat m^{\bs{\hat\beta}}}$ &  $ \hat\alpha^{\bs{\hat\beta}}_{\hat{h}^{\bs{\hat\beta}}_{CV}}$\\                  
%     \cline{3-8} 

\end{tabular}
%\caption{MISEs for the kernel estimator with a bandwidth selected by cross-validation with a Lasso estimator of the regression parameter.}
%\label{tab:MISEempLoisGL}
%\end{subtable}
%\vspace{0.4cm}
%
%
%
%\begin{subtable}[b]{0.8\linewidth}
%\begin{tabular}{|c|c|c|c|c|c|c|c|}
%\hline
%\multicolumn{2}{|c|}{\backslashbox{Dimensions}{Distributions}}  &   $\mathcal{W}(1.5,1)$  & $\mathcal{W}(0.5,2)$   &   $\mathcal{W}(3,4)$  \\
%\hline
%  \multirow{2}{*}{$n=200$} & $p=15$   &   0.072 & 0.626   & 5.26         \\
%  \cline{2-5}
%                                          & $p=200$  &  0.071   &     0.613           &      5.30          \\
% \hline
%\multirow{2}{*}{$n=500$} & $p=22$       & 0.055   &   0.401    &   5.24   \\
%\cline{2-5}
%                                        & $p=500$    &  0.059  &     0.402         &    5.25        \\                                        
%\hline
%\end{tabular}
%\caption{MISEs for the penalized contrast estimator in a histogram basis, with an adaptive Lasso estimator of the regression parameter.}
%\label{tab:MISEempLoisMS}
%\end{subtable}
%\vspace{0.5cm}

\caption{Random empirical MISE for the penalized contrast estimator in a histogram basis (first column for each distribution) and for the kernel estimator with a bandwidth selected by cross-validation (second column for each distribution), with a Lasso estimator of the regression parameter, for three different Weibull distributions of the survival times.}
\label{tab:MISEempLois}
\end{table}

First, as expected, the random $\MISE$s are smaller for a large $n$ and a small $p$. Then, we observe that the penalized contrast estimator performs better than the kernel estimator for the Weibull distributions $\mathcal{W}(0.5,2)$ and $\mathcal{W}(3,4)$. Note that the random $\MISE$s are very high for this last distribution. This can easily be explained from the fact that the baseline hazard function associated to a $\mathcal{W}(3,4)$ has the most complicated form since it increases steeply (see Figure~\ref{fig:baseline}). Lastly, for the distribution $\mathcal{W}(1.5,1)$, the random $\MISE$s are smaller in the case of the kernel estimator with a bandwidth selected by cross-validation than in the case of the penalized contrast estimator.

\section{Proofs}
\label{sec:technic2}

\subsection{Technical results}
In this section, we introduce some propositions and lemmas that are necessary to prove the theorems. Their proofs are postponed to Subsection \ref{proof:technic2}.

Let us first introduce the random norm revealed from the contrast (\ref{eq:contrast2}) and associated to the deterministic norm defined by (\ref{normdet}), and its associated scalar product: 
for $\alpha$, $\alpha_1$ and $\alpha_2$ functions in $(\mathbb{L}^{2}\cap \mathbb{L}^{\infty})([0,\tau])$ and $\bs\beta\in\mathbb{R}^p$ fixed,
\begin{eqnarray}
\label{normrand}
||\alpha||_{rand(\bs\beta)}^2\hspace{-0.2cm}&=&\hspace{-0.2cm}\dfrac{1}{n}\sumin\inttau\alpha^2(t)\e^{\bs{\beta^TZ_i}}Y_i(t)\diff t,\\
\nonumber
\langle\alpha_1,\alpha_2\rangle_{rand(\bs\beta)}&=&\dfrac{1}{n}\sumin\inttau\alpha_1(t)\alpha_2(t)\e^{\bs{\beta^TZ_i}}Y_i(t)\diff t,
\end{eqnarray}
Subsequently, to relieve the notations, we denote $||.||_{rand}:= ||.||_{rand(\bs\beta_0)} $ and the same holds for the associated scalar product. We state a key relation between $\langle.,.\rangle_{rand(\bs{\beta})}$ and $C_n(.,\bs\beta)$.
By definition, for all $m\in\mathcal{M}_n$ and $\bs\beta\in\mathbb{R}^p$,
\begin{equation}
\label{selection}
C_n(\hat\alpha^{\bs{\beta}}_{\hat m^{\bs{\beta}}},\bs\beta)+\pen(\hat m^{\bs\beta})\leq C_n(\hat\alpha^{\bs{\beta}}_m,\bs\beta)+\pen(m)\leq C_n(\alpha^{\bs{\beta_0}}_m,\bs\beta)+\pen(m),
\end{equation}
where $\hat m^{\bs\beta}=\arg\min_{m\in\mathcal{M}_n}\{C_n(\hat\alpha^{\bs\beta}_m,\bs\beta)+\pen(m)\}$.
Now, we write that 
\begin{align*}
&C_n(\hat\alpha^{\bs{\beta}}_{\hat m^{\bs{\beta}}},\bs\beta)-C_n(\alpha^{\bs{\beta_0}}_m,\bs\beta)\\
=&-\dfrac{2}{n}\sumin\inttau(\hat\alpha^{\bs{\beta}}_{\hat m^{\bs{\beta}}}-\alpha^{\bs{\beta_0}}_m)(t)\diff N_i(t)+\dfrac{1}{n}\sumin\inttau(\hat\alpha^{\bs{\beta}}_{\hat m^{\bs{\beta}}}(t)^2-\alpha^{\bs{\beta_0}}_m(t)^2)\e^{\bs{\beta^TZ_i}} Y_i(t)\diff t.
\end{align*}
Using the Doob-Meyer decomposition, we derive that
%&=-\dfrac{2}{n}\sumin\inttau(\hat\alpha^{\bs{\hat\beta}}_{\hat m^{\bs{\hat\beta}}}-\alpha^{\bs{\beta_0}}_m)(t)\Big(\diff M_i(t)+\alpha_0(t)\expbeta Y_i(t)\diff t\Big)+||\hat\alpha^{\bs{\hat\beta}}_{\hat m^{\bs{\hat\beta}}}||^2_{\rand}-||\alpha^{\bs{\beta_0}}_m||^2_{\rand}\\
\begin{align*}
&C_n(\hat\alpha^{\bs{\beta}}_{\hat m^{\bs{\beta}}},\bs\beta)-C_n(\alpha^{\bs{\beta_0}}_m,\bs\beta)\\
=&-2\langle\hat\alpha^{\bs{\beta}}_{\hat m^{\bs{\beta}}}-\alpha^{\bs{\beta_0}}_m,\alpha_0\rangle_{rand}+||\hat\alpha^{\bs{\beta}}_{\hat m^{\bs{\beta}}}||^2_{rand(\bs\beta)}-||\alpha^{\bs{\beta_0}}_m||^2_{rand(\bs\beta)}-2\nu_n(\hat\alpha^{\bs{\beta}}_{\hat m^{\bs{\beta}}}-\alpha^{\bs{\beta_0}}_m),
\end{align*}
where $\nu_n(\alpha)$ is defined by 
\begin{equation}
\label{eq:nun2}
\nu_n(\alpha)=\dfrac{1}{n}\sumin\inttau\alpha(t)\diff M_i(t).
\end{equation}
It follows that
\begin{align}
\nonumber
C_n(\hat\alpha^{\bs{\beta}}_{\hat m^{\bs{\beta}}},\bs\beta)-C_n(\alpha^{\bs{\beta_0}}_m,\bs\beta)&=||\hat\alpha^{\bs{\beta}}_{\hat m^{\bs{\beta}}}-\alpha^{\bs{\beta_0}}_m||^2_{rand(\bs\beta)}-2\nu_n(\hat\alpha^{\bs{\beta}}_{\hat m^{\bs{\beta}}}-\alpha^{\bs{\beta_0}}_m)\\
\label{eq:keyCn}
&+2\langle\hat\alpha^{\bs{\beta}}_{\hat m^{\bs{\beta}}}-\alpha^{\bs{\beta_0}}_m,\alpha^{\bs{\beta_0}}_m\rangle_{rand(\bs\beta)}-2\langle\hat\alpha^{\bs{\beta}}_{\hat m^{\bs{\beta}}}-\alpha^{\bs{\beta_0}}_m,\alpha_0\rangle_{rand}.
\end{align}

%Similar random and deterministic norms have already been introduced by  \citet{Reynaud06} to take into account the scarcity of the observations at its right-hand end. However, since \citet{Reynaud06} is interested in the estimation of the intensity function without covariates, the random norm is weighted by $Y$, when our is weighted by $\e^{\bs{\beta^TZ}}Y$.

Let us now introduce the following events :
\begin{eqnarray}
\label{set:evenements2}
\Delta_1=\left\{\alpha\in\mathcal{S}_n : \left |\dfrac{||\alpha||^2_{rand}}{||\alpha||^2_{det}}-1\right|\leq\dfrac{1}{2}\right\},\quad \mbox{and} \quad \Omega=\left\{\left|\dfrac{\hat{f}_0}{f_0}-1\right|\leq \dfrac{1}{2}\right\}
\end{eqnarray}

\begin{eqnarray}
\label{set:delta2}
\Delta_2=\left\{\alpha\in\mathcal{S}_n : \left |\dfrac{||\alpha||^2_{rand(\bs{\hat\beta})}}{||\alpha||^2_{rand}}-1\right|\leq\dfrac{1}{2}\right\}.
\end{eqnarray}
On the sets $\Delta_1$ and $\Delta_2$ we have a relation between the random $||.||_{rand}$ and the deterministic $||.||_{det}$ norms and between the random norms $||.||_{rand}$ and $||.||_{rand(\bs{\hat\beta})}$ respectively. The following proposition state a relation between the deterministic norm (\ref{normdet}) and the standard $\mathbb{L}^{2}$-norm:
\begin{proposition}[Connections between the norms]
\label{prop:connection2}
From Assumptions \ref{ass:assnot2}.\ref{ass:Zij2}-\ref{ass:f0}, we deduce the following connection between the deterministic norm and the standard $\mathbb{L}^{2}$-norm:
\begin{equation*}
\label{connec}
f_0\e^{-B|\bs{\beta_0}|_1}||\alpha||^2_2\leq ||\alpha||^2_{det}\leq \mathbb{E}[\e^{\bs{\beta_0^TZ}}]||\alpha||^2_2\leq \e^{B|\bs{\beta_0}|_1}||\alpha||^2_2.
\end{equation*}
\end{proposition}
\begin{flushleft}
The proof of this proposition is immediate using the fact that from Assumption \ref{ass:assnot2}.\ref{ass:f1}, we can rewrite the deterministic norm as 
\[
||\alpha||_{det}^2=
\inttau\int_A\alpha^2(t)e^{\bs{\beta_0^Tz}}\mathbb{E}[Y(t)|\bs{Z}=\bs{z}]f_{\bs{Z}}(z)\diff z\diff t.
\]

\end{flushleft}

\subsubsection{Results used in the proofs of Theorem \ref{th:IO2}}
Recall that for all $\bs\beta\in\mathbb{R}^p$,
\[
\hat{\mathcal{H}}^{\bs{\beta}}_m=\Bigg\{\min\Sp(\bs{G^{\bs{\beta}}_m})\geq \max\Bigg(\dfrac{\hat{f}_0\e^{-B|\bs\beta_0|_1}\e^{-B|\bs\beta_0-\bs{\beta|_1}}}{6},\dfrac{1}{\sqrt{n}}\Bigg)\Bigg\}.
\]
The following lemma ensures the existence of the estimators $\hat\alpha^{\bs{\hat\beta}}_{\hat m^{\bs{\hat\beta}}}$ on $\Delta_1\cap\Delta_2\cap\Omega$.

%\begin{lemma}
%\label{lem:incl1}
%Under Assumptions \ref{ass:assnot2}.\ref{ass:Zij2}-\ref{ass:alpha0inf2}, Assumption \ref{ass:betaball2} and Assumptions \ref{ass:model2}.\ref{ass:Dm2}-\ref{ass:emboite}, for $n\geq~ 16/(f_0\e^{-B|\bs{\beta_0}|_1})^2$, the following embedding holds:
%\[
%\Delta_1\cap\Omega\subset\hat{\mathcal{H}}\cap\Omega, \quad \mbox{ where } \hat{\mathcal{H}}:=\underset{m\in\mathcal{M}_n}{\cap}\hat{\mathcal{H}}^{\bs{\beta_0}}_m.
%\]
%\end{lemma}

\begin{lemma}
\label{lem:incl2}
Under Assumptions \ref{ass:assnot2}.\ref{ass:Zij2}-\ref{ass:alpha0inf2}, Assumptions \ref{ass:betaball2} and Assumptions \ref{ass:model2}.\ref{ass:Dm2}-\ref{ass:emboite}, for $n\geq 16/(f_0\e^{-3BR})^2$,  the following embedding holds:
\[
\Delta_1\cap\Delta_2\cap\Omega\subset\hat{\mathcal{H}}^{\bs{\hat\beta}}\cap\Omega, \quad \mbox{ where } \hat{\mathcal{H}}^{\bs{\hat\beta}}:=\underset{m\in\mathcal{M}_n}{\cap}\hat{\mathcal{H}}^{\bs{\hat\beta}}_m.
\]
\end{lemma}
\begin{flushleft}
From this lemma, for all $m\in\mathcal{M}_n$, the matrix $\bs{G^{\hat\beta}_m}$ is invertible on $\Delta_1\cap\Delta_2\cap\Omega$, and thus the estimator of $\alpha_0$ is
well defined.
Proof \ref{lem:incl2} are available in Subsection \ref{proof:incl}.

\end{flushleft}

The following proposition bounds the quadratic difference between $\hat\alpha^{\bs{\hat\beta}}_{\hat m^{\bs{\hat\beta}}}$ and $\alpha^{\bs{\beta_0}}_m$ for $m\in\mathcal{M}_n$, on the complements of 
\[
\aleph_k=\Delta_1\cap\Delta_2\cap\Omega\cap\Omega^k_H, \]
where $\Omega^k_H$, (the indice $H$ is for "Huang", since the set has already been defined by \citet{Huang2013}), is defined for $k>0$ by
\begin{align}
\label{set:OmegaH2}
\Omega^k_H=\left\{|\bs{\hat\beta}-\bs{\beta_0}|_1\leq C(s)\sqrt{\dfrac{\log (pn^k)}{n}}\right\},
\end{align}
for a constant $C(s)$ depending on the sparsity index of $\bs{\beta_0}$. From Proposition \ref{prop:predbeta2}, $\Proba(\Omega^k_H)\geq 1-cn^{-k}$ for a constant $c>0$. Now, let us state the two following propositions.

%\begin{align}
%\label{set:OmegaH2}
%\Omega^k_H=\left\{|\bs{\dot l^*_n(\bs{\beta_0})}|_{\infty} \leq \dfrac{\xi-1}{\xi+1}\Gamma_{n,k}\right\}, \quad\mbox{ with } \quad \Gamma_{n,k}=\dfrac{\xi+1}{\xi-1}B\sqrt{2\dfrac{\log(pn^k)}{n}}.
%\end{align}

%\begin{proposition}
%\label{propcomp1}
%Under Assumptions \ref{ass:assnot2}.\ref{ass:Zij2}-\ref{ass:alpha0inf2}, Assumption \ref{ass:betaball2} and Assumptions \ref{ass:model2}.\ref{ass:Dm2}-\ref{ass:emboite},
%\begin{equation}
%\mathbb{E}[||\hat\alpha^{\bs{\beta_0}}_{\hat m^{\bs{\beta_0}}}-\alpha^{\bs{\beta_0}}_m||^2_{det}\mathds{1}_{(\Delta_1\cap\Omega)^c}]\leq c_1/n,
%\end{equation}
%where $c_1$ is a constant depending on $\tau$, $\phi$, $||\alpha_0||_{\infty,\tau}$, $f_0$, $\mathbb{E}[\e^{\bs{\beta_0^TZ}}]$, $\mathbb{E}[\e^{2\bs{\beta_0^TZ}}]$, $B$, $|\bs{\beta_0}|_1$ and $\kappa_b$ a constant that comes from the B\"urkholder Inequality (see Theorem \ref{th:Burkholder}).
%\end{proposition}

%We obtain a similar result on the complement of $\aleph_k$, where $\Omega^k_H$ is defined by (\ref{set:OmegaH2}).
\begin{proposition}
\label{propcomp}
Under Assumptions \ref{ass:assnot2}.\ref{ass:Zij2}-\ref{ass:alpha0inf2}, Assumptions \ref{ass:betaball2} and Assumptions \ref{ass:model2}.\ref{ass:Dm2}-\ref{ass:emboite},
\begin{equation}
\mathbb{E}[||\hat\alpha^{\bs{\hat\beta}}_{\hat m^{\bs{\hat\beta}}}-\alpha^{\bs{\beta_0}}_m||^2_{det}\mathds{1}_{\aleph_k^c}]\leq \tilde c_1/n,
\end{equation}
where $\tilde c_1$ is a constant depending on $\tau$, $\phi$, $||\alpha_0||_{\infty,\tau}$, $f_0$, $\mathbb{E}[\e^{\bs{\beta_0^TZ}}]$, $\mathbb{E}[\e^{2\bs{\beta_0^TZ}}]$, $B$, $|\bs{\beta_0}|_1$, the sparsity index $s$ of $\bs{\beta_0}$ and $\kappa_b$ a constant that comes from the B\"urkholder Inequality (see Theorem \ref{th:Burkholder}).
\end{proposition}
\begin{flushleft}
%The proof of Propositions \ref{propcomp1} and \ref{propcomp} are similar.
We refer to Subsection \ref{proof:prop} for the proof of Proposition \ref{propcomp}. 
%A detailed proof of Proposition \ref{propcomp1} is available in \cite{CGG}.
This propositions are directly used in the proof of Theorems \ref{th:IO2} in Subsection \ref{sec:proof2}. 
%From these two results, it remains to bound $\mathbb{E}[||\hat\alpha^{\bs{\beta_0}}_{\hat m^{\bs{\beta_0}}}-\alpha^{\bs{\beta_0}}_m||^2_{det}]$ or $\mathbb{E}[||\hat\alpha^{\bs{\hat\beta}}_{\hat m^{\bs{\hat\beta}}}-\alpha^{\bs{\beta_0}}_m||^2_{det}]$ on the respective complements, to finish the proofs of Theorems \ref{th:IObetaconnu2} and \ref{th:IO2}.
%\newline

\end{flushleft}

%The following Lemma ensure that the probability of the event $\Omega^k_H$ is larger than $1-~cn^{-k}$. 
%\begin{lemma}
%\label{lem:OmegaH}
%Let consider, for $k>0$, the event $\Omega^k_H$ defined by (\ref{set:OmegaHann2}).
%Then, under Assumptions \ref{ass:Zij2} and \ref{ass:assnot2}.\ref{ass:alpha0inf2}, there exists a constant $c>0$ depending on $\tau$, $||\alpha_0||_{\infty,\tau}$ and $\mathbb{E}[\e^{\bs{\beta_0^TZ}}]$ such that 
%\[
%\mathbb{P}((\Omega^k_H)^c)\leq cn^{-k}.
%\]
%\end{lemma}
%\begin{flushleft}
%We refer to Appendix \ref{ann:Huang} for the proof of this lemma.
%\end{flushleft}
%

Usually, in model selection (see for instance \citet{Massart}), the penalty is obtained by using the so-called Talagrand's deviation inequality for the maximum of empirical processes. In the empirical process (\ref{eq:nun2}), the martingales $M_i$, $i=1,...,n$, are unbounded, 
%\begin{equation}
%\label{eq:nun2}
%\nu_n(\alpha)=\dfrac{1}{n}\sumin\inttau\alpha(t)\diff M_i(t),
%\end{equation}
Thus, we cannot directly use the Talagrand's inequality. We consider the following proposition proved in \citet{CGG}. To obtain an uniform deviation of $\nu_n(.)$, \citet{CGG} have used tools from \citet{SVG} to establish Bennett and Bernstein type inequalities and a $\mathbb{L}^{2}(det)-\mathbb{L}^{\infty}$ generic chaining type of technique (see \citet{Talagrand05generic} and \citet{Baraud10bernstein}). 

\begin{proposition}
\label{chainage}
Let $m,m'\in\mathcal{M}_n$. Define 
\begin{align}
\label{eq:BouleDet}
\mathcal{B}^{det}_{m,m'}(0,1)=\{\alpha\in S_m+S_m' : ||\alpha||_{\det}\leq1\}.
\end{align}
Under the assumptions of Theorem \ref{th:IO2}, there exists $\kappa>0$ such that for
\begin{equation}
\label{eq:pmmprime2} 
p(m,m')=\dfrac{\kappa}{K_0}(\pen(m)+\pen(m')),
\end{equation}
where the constant $K_0$ and $\pen(m)$ are defined in (\ref{eq:pen2}),
then
\[
\displaystyle\sum_{m'\in\mathcal{M}_n}\mathbb{E}\Big(\Big(\underset{\alpha\in\mathcal{B}^{det}_{m,m'}(0,1)}{\sup}\nu^2_n(\alpha)-p(m,m')\Big)_+\mathds{1}_{\Delta_1}\Big)\leq \dfrac{C_3}{n}
\]
for $n$ large enough, where $C_3$ is a constant depending on $f_0$, $\mathbb{E}[\e^{\bs{\beta_0^TZ}}]$, $B$, $|\bs{\beta_0}|_1$, $||\alpha_0||_{\infty,\tau}$ and the choice of the basis.
\end{proposition}
\begin{flushleft}
%In this proposition, the term $p(m,m')$ is the sum of two penalties up to a constant: $\forall m,m'\in\mathcal{M}_n$,
%\begin{align}
%\label{eq:pmpluspmprime2}
%p(m,m')=\dfrac{\kappa}{K_0}(\pen(m)+\pen(m')),
%\end{align}
%where $\pen(m)$ is defined by (\ref{eq:pen2}) for all $m\in\mathcal{M}_n$. 
These propositions are applied to prove Theorem \ref{th:IO2}. We admit the proof of this proposition and refer to \citet{CGG} for a detailed proof of this result. 
\end{flushleft}

We need Proposition \ref{espeta} to prove Theorem \ref{th:IO2}: the empirical centered process $\eta_n(\alpha,\alpha^{\bs{\beta_0}}_m)$, defined by 
\begin{eqnarray*}
\label{eta}
\eta_n(\alpha,\alpha^{\bs{\beta_0}}_m)=\dfrac{1}{n}\sumin\Big(U_{i}(\alpha,\alpha^{\bs{\beta_0}}_m)-\mathbb{E}[U_{i}(\alpha,\alpha^{\bs{\beta_0}}_m)]\Big),
\end{eqnarray*}
where
\[
U_{i}(\alpha,\alpha^{\bs{\beta_0}}_m)=\left(\inttau\alpha(t)\alpha^{\bs{\beta_0}}_m(t)\expbeta Y_i(t)\diff t\right)^2.
\]
appears in the proof of Theorem \ref{th:IO2}, when we control the difference between the scalar products $\langle.,.\rangle_{rand}-\hspace{-0.1cm}\langle.,.\rangle_{rand(\bs{\hat\beta})}$ (see Subsection \ref{proof:IO2}). 
Proposition \ref{espeta} allows to control this process.
\begin{proposition}
\label{espeta}
%Let consider the empirical centered process $\eta_n(\alpha,\alpha^{\bs{\beta_0}}_m)$
Let introduce the ball $\mathcal{B}^{det}_n(0,1)\subset \mathcal{S}_n$ defined by
\begin{align}
\label{eq:Bn01}
\mathcal{B}^{det}_n(0,1)=\{\alpha\in\mathcal{S}_n : ||\alpha||_{det}\leq 1\}.
\end{align}
Under Assumptions \ref{ass:assnot2}.\ref{ass:Zij2}-\ref{ass:alpha0inf2} and Assumption \ref{ass:betaball2}, we have
\[
\mathbb{E}\left[\underset{\alpha\in\mathcal{B}^{det}_n(0,1)}{\sup}\eta_n(\alpha,\alpha^{\bs{\beta_0}}_m)^2\right]\leq \dfrac{1}{n}\dfrac{\mathbb{E}[\e^{4\bs{\beta_0^TZ}}]||\alpha^{\bs{\beta_0}}_m||^4_2}{(\e^{-B|\bs{\beta_0}|_1}f_0)^2}.
\]
\end{proposition}
\begin{flushleft}
Proposition \ref{espeta} is proved in Subsection \ref{proof:espeta}.
\end{flushleft}

\subsubsection{Technical lemmas for the proofs of Proposition \ref{propcomp1} and \ref{propcomp}}
In order to prove Proposition \ref{propcomp}, we need three lemmas:

\begin{lemma}
\label{hatalpha}
Under Assumptions \ref{ass:assnot2}.\ref{ass:Zij2}-\ref{ass:alpha0inf2}, Assumptions \ref{ass:betaball2} and Assumptions \ref{ass:model2}.\ref{ass:Dm2}-\ref{ass:emboite}, we have
\[
\mathbb{E}[||\hat\alpha^{\bs{\hat\beta}}_{\hat m^{\bs{\hat\beta}}}||^4_2]\leq C_bn^4,
\]
where $C_b$ is constant depending on $||\alpha_0||_{\infty,\tau}$, $\tau$, $\mathbb{E}[\e^{\bs{\beta_0^TZ}}]$ and $\mathbb{E}[\e^{2\bs{\beta_0^TZ}}]$, $\kappa_b$, the constant of the B\"urkholder Inequality (see Theorem \ref{th:Burkholder}) and on the choice of the basis.
\end{lemma}

\begin{lemma}
\label{Delta1}
Under Assumptions \ref{ass:assnot2}.\ref{ass:Zij2}-\ref{ass:alpha0inf2} and Assumptions \ref{ass:model2}.\ref{ass:Dm2}-\ref{ass:emboite}, we have
\[
\mathbb{P}(\Delta^c_1)\leq \dfrac{C^{(\Delta_1)}_k}{n^k}, \quad \forall k\geq1,
\]
where $C^{(\Delta_1)}_k$ is a constant depending on $f_0$, $B$ and $|\bs{\beta_0}|_1$.
\end{lemma}

\begin{lemma}
\label{Delta2}
Under Assumptions \ref{ass:assnot2}.\ref{ass:Zij2}-\ref{ass:alpha0inf2}, Assumptions \ref{ass:betaball2} and Assumption \ref{ass:Relationpn}, we have for $n$ large enough,
\[
\mathbb{P}(\Delta^c_2)\leq \dfrac{C^{(\Delta_2)}_k}{n^{k}}, \quad \forall k\geq 1,
\]
where the constant $C^{(\Delta_2)}_k$ depends on $\tau$, $||\alpha_0||_{\infty,\tau}$ and $\mathbb{E}[\e^{\bs{\beta_0^TZ}}]$.
\end{lemma}

These three lemmas are required to prove Proposition \ref{propcomp}. There are proved in Subsection \ref{proof:technic2}.

\subsubsection{A classical inequality: the B\"urkholder Inequality}

The last technical result is a B\"urkholder Inequality that gives a norm relation between a martingale and its optional process. We refer to \citet{Liptser89} p.75, for the proof of this result.
\begin{theorem}[B\"urkholder Inequality] 
\label{th:Burkholder}
If $M=(M_t,\mathcal{F}_t)_{t\geq 0}$ is a martingale, then there are universal constants $\gamma_b$ and $\kappa_b$ (independent of $M$) such that for every $t\geq 0$
\[
\gamma_b||\sqrt{[M]_t}||_2\leq||M_t||_2\leq \kappa_b||\sqrt{[M]_t}||_2,
\]
where $[M]_t$ is the quadratic variation of $M_t$. 
%The constants $\gamma_b$ and $\kappa_b$ can be taken to have the values 
%\[
%\gamma_b=[18b^{3/2}/(b-1)]^{-1}, \quad \kappa_b=18b^{3/2}/(b-1)^{1/2}.
%\]
\end{theorem}
\begin{flushleft}
This theorem is used to prove Lemma \ref{hatalpha} and in the oracle inequalities of Theorem \ref{th:IO2}, the constants depend on  $\kappa_b$.
\end{flushleft}

%%% SECTION 4 --- PROOFS OF THE MAIN THEOREMS ---
\subsection{Proofs of the main theorems}
\label{sec:proof2}

\subsubsection{Proof of Theorem \ref{th:IO2}}
\label{proof:IO2}
In the following, we consider the sets $\Delta_1$, $\Delta_2$ and $\Omega$ defined by (\ref{set:evenements2}) and (\ref{set:delta2}) and the set $\Omega^k_H$ defined by (\ref{set:OmegaH2}). For sake of simplicity in the notations, we denote $\aleph_k$ the intersection between the four sets: $\aleph_k=\Delta_1\cap\Delta_2\cap\Omega\cap\Omega^k_H$. 
We have the following decomposition:
\begin{align*}
\mathbb{E}[||\hat\alpha^{\bs{\hat\beta}}_{\hat m^{\bs{\hat\beta}}}-\alpha_0||^2_{det}]\leq 2||\alpha_0-\alpha^{\bs{\beta_0}}_m||^2_{det}+2\mathbb{E}[||\hat\alpha^{\bs{\hat\beta}}_{\hat m^{\bs{\hat\beta}}}-\alpha^{\bs{\beta_0}}_m||^2_{det}\mathds{1}_{\aleph_k}]
+2\mathbb{E}[||\hat\alpha^{\bs{\hat\beta}}_{\hat m^{\bs{\hat\beta}}}-\alpha^{\bs{\beta_0}}_m||^2_{det}\mathds{1}_{\aleph_k^c}].
\end{align*}
The first term is the usual bias term.
From Proposition \ref{propcomp}, we deduce that the last term is bounded by $\tilde c_1/n$.
We now focus on the term $\mathbb{E}[||\hat\alpha^{\bs{\hat\beta}}_{\hat m^{\bs{\hat\beta}}}-\alpha^{\bs{\beta_0}}_m||^2_{det}\mathds{1}_{\aleph_k}]$. 
From Lemma \ref{lem:incl2}, for all $m\in\mathcal{M}_n$, the matrices $\bs{G^{\bs{\hat\beta}}_m}$ are invertible on $\Delta_1\cap \Delta_2 \cap \Omega\cap\Omega^k_H$ as soon as $n\geq 16/(f_0\e^{-3BR})^2$ and thus the estimator $\hat\alpha^{\bs{\hat\beta}}_{\hat m^{\bs{\hat\beta}}}$ of $\alpha_0$ is well defined. 
From (\ref{selection}) and (\ref{eq:keyCn}), with $\bs\beta=\bs{\hat\beta}$, we have for all $m\in\mathcal{M}_n$,
%$||\hat\alpha^{\bs{\hat\beta}}_{\hat m^{\bs{\hat\beta}}}-\alpha^{\bs{\beta_0}}_m||^2_{\rand}$ is less than
\begin{align*}
%&\leq 2\nu_n(\hat\alpha^{\bs{\hat\beta}}_{\hat m^{\bs{\hat\beta}}}-\alpha^{\bs{\beta_0}}_m)+\pen(m)-\pen(\hat m^{\hatbeta})\\
%&\hspace{2.65cm}2\langle\hat\alpha^{\bs{\hat\beta}}_{\hat m^{\bs{\hat\beta}}}-\alpha^{\bs{\beta_0}}_m,\alpha_0\rangle_{rand}-2\langle\hat\alpha^{\bs{\hat\beta}}_{\hat m^{\bs{\hat\beta}}}-\alpha^{\bs{\beta_0}}_m,\alpha^{\bs{\beta_0}}_m\rangle_{\rand}\\
||\hat\alpha^{\bs{\hat\beta}}_{\hat m^{\bs{\hat\beta}}}&-\alpha^{\bs{\beta_0}}_m||^2_{\rand}\leq  2\nu_n(\hat\alpha^{\bs{\hat\beta}}_{\hat m^{\bs{\hat\beta}}}-\alpha^{\bs{\beta_0}}_m)+2\langle\hat\alpha^{\bs{\hat\beta}}_{\hat m^{\bs{\hat\beta}}}-\alpha^{\bs{\beta_0}}_m,\alpha_0-\alpha^{\bs{\beta_0}}_m\rangle_{rand}\\
 &+\pen(m)-\pen(\hat m^{\hatbeta})+2\langle\hat\alpha^{\bs{\hat\beta}}_{\hat m^{\bs{\hat\beta}}}-\alpha^{\bs{\beta_0}}_m,\alpha^{\bs{\beta_0}}_m\rangle_{rand}-2\langle\hat\alpha^{\bs{\hat\beta}}_{\hat m^{\bs{\hat\beta}}}-\alpha^{\bs{\beta_0}}_m,\alpha^{\bs{\beta_0}}_m\rangle_{\rand},
\end{align*}
where the empirical process $\nu_n(.)$ is defined by Equation (\ref{eq:nun2}) and the random norm by (\ref{normrand}). For $\mathcal{B}^{det}_{m,m'}(0,1)$ defined by (\ref{eq:BouleDet}), using the classical inequality $2xy\leq ~bx^2+y^2/b$ with $b>0$, we obtain
\begin{align*}
||\hat\alpha^{\bs{\hat\beta}}_{\hat m^{\bs{\hat\beta}}}-\alpha^{\bs{\beta_0}}_m||^2_{\rand}&
%\leq \dfrac{1}{16}||\hat\alpha^{\bs{\hat\beta}}_{\hat m^{\bs{\hat\beta}}}-\alpha^{\bs{\beta_0}}_m||^2_{rand}+16||\alpha_0-\alpha^{\bs{\beta_0}}_m||^2_{rand}+\pen(m)-\pen(\hat m^{\hatbeta})\\
%&\hspace{0.7cm}+2\nu_n(\hat\alpha^{\bs{\hat\beta}}_{\hat m^{\bs{\hat\beta}}}-\alpha^{\bs{\beta_0}}_m)+2\Big(\langle\hat\alpha^{\bs{\hat\beta}}_{\hat m^{\bs{\hat\beta}}}-\alpha^{\bs{\beta_0}}_m,\alpha^{\bs{\beta_0}}_m\rangle_{rand}-\langle\hat\alpha^{\bs{\hat\beta}}_{\hat m^{\bs{\hat\beta}}}-\alpha^{\bs{\beta_0}}_m,\alpha^{\bs{\beta_0}}_m\rangle_{\rand}\Big)\\
\leq  \dfrac{1}{16}||\hat\alpha^{\bs{\hat\beta}}_{\hat m^{\bs{\hat\beta}}}-\alpha^{\bs{\beta_0}}_m||^2_{rand}+16||\alpha_0-\alpha^{\bs{\beta_0}}_m||^2_{rand}+\pen(m)-\pen(\hat m^{\hatbeta})\\
&+\dfrac{1}{16}||\hat\alpha^{\bs{\hat\beta}}_{\hat m^{\bs{\hat\beta}}}-\alpha^{\bs{\beta_0}}_m||^2_{det}+16\underset{\alpha\in\mathcal{B}^{det}_{m,\hat m^{\hatbeta}}(0,1)}{\sup}\nu^2_n(\alpha)\\
&+2\Big(\langle\hat\alpha^{\bs{\hat\beta}}_{\hat m^{\bs{\hat\beta}}}-\alpha^{\bs{\beta_0}}_m,\alpha^{\bs{\beta_0}}_m\rangle_{rand}-\langle\hat\alpha^{\bs{\hat\beta}}_{\hat m^{\bs{\hat\beta}}}-\alpha^{\bs{\beta_0}}_m,\alpha^{\bs{\beta_0}}_m\rangle_{\rand}\Big).
\end{align*}
%\[
%\mathbb{E}[||\hat\alpha^{\bs{\beta_0}}_{\hat m^{\bs{\beta_0}}}-\alpha_0||^2_{det}]\leq 2||\alpha_0-\alpha^{\bs{\beta_0}}_m||^2_{det}+2\mathbb{E}[||\hat\alpha^{\bs{\beta_0}}_{\hat m^{\bs{\beta_0}}}-\alpha^{\bs{\beta_0}}_m||^2_{det}\mathds{1}_{(\Delta_1\cap\Delta_2\cap\Omega)}]+2\mathbb{E}[||\hat\alpha^{\bs{\beta_0}}_{\hat m^{\bs{\beta_0}}}-\alpha^{\bs{\beta_0}}_m||^2_{det}\mathds{1}_{(\Delta_1\cap\Delta_2\cap\Omega)^c}]
%\]
Consequently, using the relations between the random norms $||.||_{rand(\bs{\hatbeta})}$ and $||.||_{rand}$ and between the random norm $||.||_{rand}$ and the deterministic norm $||.||_{det}$ on $\aleph_k$, we obtain
\begin{align*}
\dfrac{1}{4}||\hat\alpha^{\bs{\hat\beta}}_{\hat m^{\bs{\hat\beta}}}-\alpha^{\bs{\beta_0}}_m||^2_{det}&\leq \dfrac{3}{32}||\hat\alpha^{\bs{\hat\beta}}_{\hat m^{\bs{\hat\beta}}}-\alpha^{\bs{\beta_0}}_m||^2_{det}+16||\alpha_0-\alpha^{\bs{\beta_0}}_m||^2_{rand}+\pen(m)-\pen(\hat m^{\hatbeta})\\
&\hspace{0.7cm}+\dfrac{1}{16}||\hat\alpha^{\bs{\hat\beta}}_{\hat m^{\bs{\hat\beta}}}-\alpha^{\bs{\beta_0}}_m||^2_{det}+16\underset{\alpha\in\mathcal{B}^{det}_{m,\hat m^{\hatbeta}}(0,1)}{\sup}\nu^2_n(\alpha)\\
&\hspace{0.7cm}+2\Big(\langle\hat\alpha^{\bs{\hat\beta}}_{\hat m^{\bs{\hat\beta}}}-\alpha^{\bs{\beta_0}}_m,\alpha^{\bs{\beta_0}}_m\rangle_{rand}-\langle\hat\alpha^{\bs{\hat\beta}}_{\hat m^{\bs{\hat\beta}}}-\alpha^{\bs{\beta_0}}_m,\alpha^{\bs{\beta_0}}_m\rangle_{\rand}\Big),
\end{align*}
also be rewritten for $p(m,m')$ defined by (\ref{eq:pmmprime2}) for all $m'\in\mathcal{M}_n$, as
\begin{align*}
&\dfrac{3}{32}\mathbb{E}\Big[||\hat\alpha^{\bs{\hat\beta}}_{\hat m^{\bs{\hat\beta}}}-\alpha^{\bs{\beta_0}}_m||^2_{det}\mathds{1}_{\aleph_k}\Big]\leq 16||\alpha_0-\alpha^{\bs{\beta_0}}_m||^2_{det} +16p(m,\hat m^{\bs{\hat\beta}})\\
&+\pen(m)-\pen(\hat m^{\bs{\hat\beta}})+16\displaystyle\sum_{m'\in\mathcal{M}_n}\mathbb{E}\Bigg(\Bigg(\underset{\alpha\in\mathcal{B}^{det}_{m,m'}(0,1)}{\sup}\nu^2_n(\alpha)-p(m,m')\Bigg)_+\mathds{1}_{\aleph_k}\Bigg)\\
&+ 2 \mathbb{E}\Big[\Big( \langle\hat\alpha^{\bs{\hat\beta}}_{\hat m^{\bs{\hat\beta}}}-\alpha^{\bs{\beta_0}}_m,\alpha^{\bs{\beta_0}}_m\rangle_{rand}-\langle\hat\alpha^{\bs{\hat\beta}}_{\hat m^{\bs{\hat\beta}}}-\alpha^{\bs{\beta_0}}_m,\alpha^{\bs{\beta_0}}_m\rangle_{\rand}\Big)\mathds{1}_{\aleph_k}\Big].
\end{align*}
We fix $K_0\geq 16\kappa$ such that $16p(m,m')\leq \pen(m)+\pen(m')$, for all $m,m'$ in $\mathcal{M}_n$, so that
\begin{align*}
\dfrac{3}{32}\mathbb{E}\Big[||\hat\alpha^{\bs{\hat\beta}}_{\hat m^{\bs{\hat\beta}}}-\alpha^{\bs{\beta_0}}_m||^2_{det}&\mathds{1}_{\aleph_k}\Big]\leq 16||\alpha_0-\alpha^{\bs{\beta_0}}_m||^2_{det}+2\pen(m)\\
&+16\displaystyle\sum_{m'\in\mathcal{M}_n}\mathbb{E}\Big(\Big(\underset{\alpha\in\mathcal{B}^{det}_{m,m'}(0,1)}{\sup}\nu^2_n(\alpha)-p(m,m')\Big)_+\mathds{1}_{\aleph_k}\Big)\\
&+2 \mathbb{E}\Big[\Big( \langle\hat\alpha^{\bs{\hat\beta}}_{\hat m^{\bs{\hat\beta}}}-\alpha^{\bs{\beta_0}}_m,\alpha^{\bs{\beta_0}}_m\rangle_{rand}-\langle\hat\alpha^{\bs{\hat\beta}}_{\hat m^{\bs{\hat\beta}}}-\alpha^{\bs{\beta_0}}_m,\alpha^{\bs{\beta_0}}_m\rangle_{\rand}\Big)\mathds{1}_{\aleph_k}\Big],
\end{align*}
that is
\begin{align}
\label{eq:io} 
\dfrac{3}{32}\mathbb{E}[||\hat\alpha^{\bs{\hat\beta}}_{\hat m^{\bs{\hat\beta}}}-\alpha^{\bs{\beta_0}}_m||^2_{det}\mathds{1}_{\aleph_k}]\leq 16||\alpha_0-\alpha^{\bs{\beta_0}}_m||^2_{det}+2\pen(m)+A(m)+\mathbb{E}[B(m,\hat m^{\hatbeta})\mathds{1}_{\aleph_k}]
%\\&+4p(m,\hat m(\hatbeta))-\pen(\hat m(\hatbeta)),
\end{align}
where
\begin{align}
\label{nu}
A(m)&=16\displaystyle\sum_{m'\in\mathcal{M}_n}\mathbb{E}\Bigg(\Bigg(\underset{\alpha\in\mathcal{B}^{det}_{m,m'}(0,1)}{\sup}\nu^2_n(\alpha)-p(m,m')\Bigg)_+\mathds{1}_{\aleph_k}\Bigg),\\
B(m,\hat m^{\hatbeta})
&=2\Big(\langle\hat\alpha^{\bs{\hat\beta}}_{\hat m^{\bs{\hat\beta}}}-\alpha^{\bs{\beta_0}}_m,\alpha^{\bs{\beta_0}}_m\rangle_{rand}-\langle\hat\alpha^{\bs{\hat\beta}}_{\hat m^{\bs{\hat\beta}}}-\alpha^{\bs{\beta_0}}_m,\alpha^{\bs{\beta_0}}_m\rangle_{\rand}\Big). \label{prodscal}
\end{align}
It remains to study the terms $A(m)$ and $B(m,\hat m^{\hatbeta})$.
%We have now to control two terms :
%\begin{enumerate}
%\item \label{nu}$\underset{\alpha\in\mathcal{B}^{det}_{m,\hat m}(0,1)}{\sup}\nu^2_n(\alpha)$
%\item \label{prodscal} $2(<\hat\alpha^{\bs{\beta_0}}_{\hat m^{\bs{\beta_0}}}-\alpha^{\bs{\beta_0}}_m,\alpha^{\bs{\beta_0}}_m>_{rand}-<\hat\alpha^{\bs{\beta_0}}_{\hat m^{\bs{\beta_0}}}-\alpha^{\bs{\beta_0}}_m,\alpha^{\bs{\beta_0}}_m>_{\rand})$
%\end{enumerate}

\paragraph{Study of (\ref{nu}).}
According to Proposition \ref{chainage}, for $n$ large enough  
\[
\displaystyle\sum_{m'\in\mathcal{M}_n}\mathbb{E}
\Bigg(\Bigg(\underset{\alpha\in\mathcal{B}^{det}_{m,m'}(0,1)}{\sup}\nu^2_n(\alpha)-p(m,m')\Bigg)_+\mathds{1}_{\aleph_k}\Bigg)\leq \dfrac{C_3}{n},
\]
where $p(m,m')$ is defined by (\ref{eq:pmmprime2}) and $C_3$ is a constant depending on $f_0$, $|\bs{\beta_0}|_1$, $B$, $\E[\e^{\bs{\beta_0^TZ}}]$, $||\alpha_0||_{\infty,\tau}$ and the choice of the basis.
Hence, for $C'_3=16C_3$, we conclude that
\begin{eqnarray}
\label{Resnu}
A(m) \leq \dfrac{C'_3}{n}.
\end{eqnarray}
%where $C'_3=16C_3$.

\paragraph{Study of (\ref{prodscal}).}
Using again the classical inequality $2xy\leq bx^2+y^2/b$ with $b>0$, we obtain
\begin{align}
\nonumber
\langle&\hat\alpha^{\bs{\hat\beta}}_{\hat m^{\bs{\hat\beta}}}-\alpha^{\bs{\beta_0}}_m,\alpha^{\bs{\beta_0}}_m\rangle_{rand}-\langle\hat\alpha^{\bs{\hat\beta}}_{\hat m^{\bs{\hat\beta}}}-\alpha^{\bs{\beta_0}}_m,\alpha^{\bs{\beta_0}}_m\rangle_{\rand}
%&=\dfrac{1}{n}\sumin\inttau(\hat\alpha^{\bs{\hat\beta}}_{\hat m^{\bs{\hat\beta}}}-\alpha^{\bs{\beta_0}}_m)(t)\alpha^{\bs{\beta_0}}_m(t)(\expbeta-\exphatbeta)Y_i(t)\diff t\\
\nonumber
\leq \dfrac{1}{32}||\hat\alpha^{\bs{\hat\beta}}_{\hat m^{\bs{\hat\beta}}}-\alpha^{\bs{\beta_0}}_m||^2_{det}\\
\label{eq:B}
&+32\underset{\alpha\in\mathcal{B}^{det}_{m,\hat m^{\hat\beta}}(0,1)}{\sup}\Big( \dfrac{1}{n}\sumin\inttau\alpha(t)\alpha^{\bs{\beta_0}}_m(t)(\expbeta-\exphatbeta)Y_i(t)\diff t\Big)^2.
%&\leq \dfrac{1}{32}||\hat\alpha^{\bs{\beta_0}}_{\hat m^{\bs{\beta_0}}}-\alpha^{\bs{\beta_0}}_m||^2_{det}+32\underset{\alpha\in\mathcal{B}^{det}_{m,\hat m}(0,1)}{\sup}\Big( \dfrac{1}{n}\sumin\inttau\alpha(t)\alpha^{\bs{\beta_0}}_m(t)(\expbeta-\exphatbeta)Y_i(t)\diff t\Big)^2,
\end{align}
%where $\mathcal{B}^{det}_{m,m'}(0,1)=\{\alpha\in S_m+S_m' : ||\alpha||_{\det}\leq 1\}$.
Now, from Assumption \ref{ass:model2}.\ref{ass:emboite} and by definition (\ref{eq:Bn01}) of $\mathcal{B}^{det}_n(0,1)$, we write that 
\[
\underset{\alpha\in\mathcal{B}^{det}_{m,\hat m^{\hat\beta}}(0,1)}{\sup}\Bigg(\dfrac{1}{n}\sumin\inttau\alpha(t)\alpha^{\bs{\beta_0}}_m(t)(\expbeta-\exphatbeta)Y_i(t)\diff t\Bigg)^2
\]
is less than
\begin{align*}
%\underset{\alpha\in\mathcal{B}^{det}_{m,\hat m}(0,1)}{\sup}\Big(\dfrac{1}{n}\sumin\inttau\alpha(t)\alpha^{\bs{\beta_0}}_m(t)(\expbeta-\exphatbeta)&Y_i(t)\diff t\Big)^2\leq
%\underset{\alpha\in\mathcal{B}^{det}_n(0,1)}{\sup}&\Big(\dfrac{1}{n}\sumin\inttau\alpha(t)\alpha^{\bs{\beta_0}}_m(t)(\expbeta-\exphatbeta)Y_i(t)\diff t\Big)^2\\
\underset{\alpha\in\mathcal{B}^{det}_n(0,1)}{\sup}\Bigg(\dfrac{1}{n}\sumin\inttau\alpha(t)\alpha^{\bs{\beta_0}}_m(t)\expbeta(1-\e^{\bs{\hat\beta^TZ_i}-\bs{\beta_0^TZ_i}}) Y_i(t)\diff t\Bigg)^2.\\
%&\leq |1-\expnorm|^2\underset{\alpha\in\mathcal{B}^{det}_n(0,1)}{\sup}\dfrac{1}{n}\sumin\Big(\inttau\alpha(t)\alpha^{\bs{\beta_0}}_m(t)\expbeta Y_i(t)\diff t\Big)^2,
%\leq \displaystyle\sum_{j\in J_m}\Big(\dfrac{1}{n}\sumin\inttau\varphi_j(t)\alpha^{\bs{\beta_0}}_m(t)\expbeta(1-\e^{\bs{\hat\beta^TZ_i}-\bs{\beta_0^TZ_i}})Y_i(t)\diff t\Big)^2\\
%&\leq (1-\expBnorm)^2\displaystyle\sum_{j\in J_m}\Big(\dfrac{1}{n}\sumin\sqrt{\inttau\varphi^2_j(t)\diff t}\sqrt{\inttau\alpha^2_m(t)\e^{2\bs{\beta_0^TZ_i}}Y^2_i(t)\diff t}\Big)^2\\
%&\leq (1-\expBnorm)^2\displaystyle\sum_{j\in J_m}\dfrac{1}{n}\sumin\inttau\alpha^2_m(t)\e^{2\bs{\beta_0^TZ_i}}Y^2_i(t)\diff t\\
%&\leq (1-\expBnorm)^2\e^{B||\bs{\beta_0}||_1}\displaystyle\sum_{j\in J_m}||\alpha^{\bs{\beta_0}}_m||^2_{rand}\\
%&\leq \e^{B|\bs\beta_0|_1}|m|||\alpha^{\bs{\beta_0}}_m||^2_{rand}(1-\expBnorm)^2
\end{align*}
%where $\mathcal{B}^{det}_n(0,1)$ is defined by (\ref{eq:Bn01}).
%\begin{align}
%\label{eq:Bn01}
%\mathcal{B}^{det}_n(0,1)=\{\alpha\in\mathcal{S}_n : ||\alpha||_{det}\leq 1\}.
%\end{align}
We have
\begin{align*}
\Bigg|\dfrac{1}{n}\sumin&\inttau\alpha(t)\alpha^{\bs{\beta_0}}_m(t)\expbeta(1-\e^{\bs{\hat\beta^TZ_i}-\bs{\beta_0^TZ_i}}) Y_i(t)\diff t\Bigg|\\
&\leq \dfrac{1}{n}\sumin\Big|1- \e^{\bs{\hat\beta^TZ_i}-\bs{\beta_0^TZ_i}}\Big|\Bigg|\inttau\alpha(t)\alpha^{\bs{\beta_0}}_m(t)\expbeta Y_i(t)\diff t\Bigg|.
\end{align*}
Using the fact that $|\e^x-\e^y|\leq|x-y|\e^{x\vee y}$ for all $(x,y)\in\mathbb{R}^2$ and applying Assumptions \ref{ass:assnot2}.\ref{ass:Zij2} and Assumptions \ref{ass:betaball2}, we obtain that
\begin{align*}
\Bigg|\dfrac{1}{n}\sumin&\inttau\alpha(t)\alpha^{\bs{\beta_0}}_m(t)\expbeta(1-\e^{\bs{\hat\beta^TZ_i}-\bs{\beta_0^TZ_i}}) Y_i(t)\diff t\Bigg|\\
&\leq \dfrac{1}{n}\sumin|\bs{\hat\beta^TZ_i}-\bs{\beta_0^TZ_i}|\e^{|\bs{\hat\beta^TZ_i}-\bs{\beta_0^TZ_i}|}\Bigg|\inttau\alpha(t)\alpha^{\bs{\beta_0}}_m(t)\expbeta Y_i(t)\diff t\Bigg|\\
&\leq B\e^{2BR}|\bs{\hat\beta}-\bs{\beta_0}|_1\Bigg|\inttau\alpha(t)\alpha^{\bs{\beta_0}}_m(t)\expbeta Y_i(t)\diff t\Bigg|.
\end{align*}
Now, write
%In order to simplify the notations in the calculations, we introduce the processes $\eta_n(\alpha,\alpha^{\bs{\beta_0}}_m)$ and $D_n(\alpha^{\bs{\beta_0}}_m)$ such that
\begin{align}
\nonumber
&\underset{\alpha\in\mathcal{B}^{det}_{n}(0,1)}{\sup}\Bigg(\dfrac{1}{n}\sumin\inttau\alpha(t)\alpha^{\bs{\beta_0}}_m(t)\expbeta(1-\e^{\bs{\hat\beta^TZ_i}-\bs{\beta_0^TZ_i}}) Y_i(t)\diff t\Bigg)^2\\
\nonumber
\leq & B^2\e^{4BR}|\bs{\hat\beta}-\bs{\beta_0}|^2_1\underset{\alpha\in\mathcal{B}^{det}_{n}(0,1)}{\sup}\dfrac{1}{n}\sumin\Bigg(\inttau\alpha(t)\alpha^{\bs{\beta_0}}_m(t)\expbeta Y_i(t)\diff t\Bigg)^2\\
\label{eq:eta}
\leq &B^2\e^{4BR}|\bs{\hat\beta}-\bs{\beta_0}|^2_1\underset{\alpha\in\mathcal{B}^{det}_{n}(0,1)}{\sup}\{\eta_n(\alpha,\alpha^{\bs{\beta_0}}_m)+D_n(\alpha,\alpha^{\bs{\beta_0}}_m)\}
%\mathbb{E}\Bigg[\Big|1-\expnorm\Big|^2\Bigg\{
%\underset{\alpha\in\mathcal{B}^{det}_n(0,1)}{\sup} \eta_n(\alpha,\alpha^{\bs{\beta_0}}_m)+ D_n(\alpha^{\bs{\beta_0}}_m)\Bigg\}\Bigg] 
%\dfrac{1}{n}\sumin\Big[\Big(\inttau\alpha(t)\alpha^{\bs{\beta_0}}_m(t)\expbeta Y_i(t)\diff t\Big)^2-\mathbb{E}\Big[\Big(\inttau\alpha(t)\alpha^{\bs{\beta_0}}_m(t)\expbeta Y_i(t)\diff t\Big)^2\Big]\Big]\\
\end{align}
where $\eta_n(\alpha,\alpha^{\bs{\beta_0}}_m)$ is defined by
\begin{equation*}
\label{eta}
\eta_n(\alpha,\alpha^{\bs{\beta_0}}_m)\hspace{-0.05cm}=\hspace{-0.05cm}\dfrac{1}{n}\sumin\Bigg[\Bigg(\inttau\alpha(t)\alpha^{\bs{\beta_0}}_m(t)\expbeta Y_i(t)\diff t\Bigg)^2\hspace{-0.2cm}-\mathbb{E}\Bigg[\Bigg(\inttau\alpha(t)\alpha^{\bs{\beta_0}}_m(t)\expbeta Y_i(t)\diff t\Bigg)^2\Bigg]\Bigg],
\end{equation*}
%with 
%\[
%U_{n,i}=\left(\inttau\alpha(t)\alpha^{\bs{\beta_0}}_m(t)\expbeta Y_i(t)\diff t\right)^2,
%\]
and 
\begin{eqnarray*}
\label{second}
D_n(\alpha,\alpha^{\bs{\beta_0}}_m)=\mathbb{E}\left[\left(\inttau\alpha(t)\alpha^{\bs{\beta_0}}_m(t)\e^{\bs{\beta_0^TZ}} Y(t)\diff t\right)^2\right].
\end{eqnarray*}
We first claim that the term $\sup_{\alpha\in\mathcal{B}^{det}_{n}(0,1)}\{D_n(\alpha,\alpha^{\bs{\beta_0}}_m)\}$ is bounded, by using that from the Cauchy-Schwarz Inequality, 
\begin{equation*}
\label{secondBound}
\underset{\alpha\in\mathcal{B}^{det}_n(0,1)}{\sup}\mathbb{E}\left[\left(\inttau\alpha(t)\alpha^{\bs{\beta_0}}_m(t)\e^{\bs{\beta_0^TZ}} Y(t)\diff t\right)^2\right]\leq ||\alpha^{\bs{\beta_0}}_m||^2_{det}.
\end{equation*}
Thus, gathering bounds (\ref{eq:B}) and (\ref{eq:eta}, we obtain that 
\[
B(m,\hat m^{\bs{\hat\beta}})\leq \dfrac{1}{16}||\hat\alpha^{\bs{\hat\beta}}_{\hat m^{\bs{\hat\beta}}}-\alpha^{\bs{\beta_0}}_m||^2_{det}+64\Bigg[B^2\e^{4BR}|\bs{\hat\beta}-\bs{\beta_0}|^2_1\Big(\underset{\alpha\in\mathcal{B}^{det}_{n}(0,1)}{\sup}\{\eta_n(\alpha,\alpha^{\bs{\beta_0}}_m)\}+||\alpha^{\bs{\beta_0}}_m||^2_{det}\Big)\Bigg].
\]
%Now, by definition (\ref{def:alpham}) of $\alpha^{\bs{\beta_0}}_m$ and from Proposition \ref{prop:connection2}, $||\alpha^{\bs{\beta_0}}_m||^2_{det}\leq||\alpha_0||^2_{det}\leq \mathbb{E}[\e^{\bs{\beta_0^TZ}}]||\alpha_0||^2_2$, and we obtain
So, taking the expectation and applying Proposition \ref{espeta} to control 
\[
\mathbb{E}[{\sup}_{\alpha\in\mathcal{B}^{det}_n(0,1)}(\eta_n(\alpha,\alpha^{\bs{\beta_0}}_m))^2],
\]
 we get
\begin{align}
\nonumber
&\E[B(m,\hat m^{\bs{\hat\beta}})\mathds{1}_{\aleph_k}]\leq \dfrac{1}{16}\E[||\hat\alpha^{\bs{\hat\beta}}_{\hat m^{\bs{\hat\beta}}}-\alpha^{\bs{\beta_0}}_m||^2_{det}\mathds{1}_{\aleph_k}]\\
\label{Resprodscal}
+&64B^2\e^{4BR}\Bigg\{\E^{1/2}[|\bs{\hat\beta}-\bs{\beta_0}|^4_1\mathds{1}_{\aleph_k}]\E^{1/2}\Bigg[\underset{\alpha\in\mathcal{B}^{det}_{n}(0,1)}{\sup}\{\eta^2_n(\alpha,\alpha^{\bs{\beta_0}}_m)\}\Bigg]\hspace{-0.1cm}+\hspace{-0.1cm}||\alpha^{\bs{\beta_0}}_m||^2_{det}\E[|\bs{\hat\beta}-\bs{\beta_0}|^2_1\mathds{1}_{\aleph_k}]\hspace{-0.1cm}\Bigg\}.
\end{align}

Finally, combining (\ref{eq:io}), (\ref{Resnu}) and (\ref{Resprodscal}) we conclude that
\begin{align*}
\dfrac{1}{16}\mathbb{E}[||\hat\alpha^{\bs{\hat\beta}}_{\hat m^{\bs{\hat\beta}}}-\alpha^{\bs{\beta_0}}_m||^2_{det}\mathds{1}_{\aleph_k}]&\leq 16||\alpha_0-\alpha^{\bs{\beta_0}}_m||^2_{det}+2\pen(m)+\dfrac{C'_3}{n}\\
&+64B^2\e^{4BR} ||\alpha^{\bs{\beta_0}}_m||^2_{det}\E[|\bs{\hat\beta}-\bs{\beta_0}|^2_1\mathds{1}_{\aleph_k}]\\
&+64B^2\e^{4BR}\E^{1/2}[|\bs{\hat\beta}-\bs{\beta_0}|^4_1\mathds{1}_{\aleph_k}]\dfrac{\mathbb{E}^{1/2}[\e^{4\bs{\beta_0^TZ}}]||\alpha^{\bs{\beta_0}}_m||^2_2}{\e^{-B|\bs{\beta_0}|_1}f_0}\dfrac{1}{\sqrt{n}}.
\end{align*}
On $\Omega\cap\Omega^k_H$, using that, from definition (\ref{def:alpham}) and Proposition \ref{prop:connection2}, $||\alpha^{\bs{\beta_0}}_m||^2_{det}\leq 2||\alpha_0||_{det}\leq\E[\e^{\bs{\beta_0}^TZ}]\tau||\alpha_0||_{\infty,\tau}$, we have
\begin{align*}
64B^2\e^{4BR} ||\alpha^{\bs{\beta_0}}_m||^2_{det}\E[|\bs{\hat\beta}-\bs{\beta_0}|^2_1\mathds{1}_{\aleph_k}]\leq C(s,B,R,\E[\e^{\bs{\beta_0^TZ}}], ||\alpha_0||_{\infty,\tau},\tau)\dfrac{\log (pn^k)}{n},
\end{align*}
and that
\begin{align*}
64B^2\e^{4BR}&\E^{1/2}[|\bs{\hat\beta}-\bs{\beta_0}|^4_1\mathds{1}_{\aleph_k}]\dfrac{\mathbb{E}^{1/2}[\e^{4\bs{\beta_0^TZ}}]||\alpha^{\bs{\beta_0}}_m||^2_2}{\e^{-B|\bs{\beta_0}|_1}f_0}\dfrac{1}{\sqrt{n}}\\
&\leq \tilde{C}(s,B,|\bs{\beta_0}|_1,R,\mathbb{E}[\e^{\bs{\beta_0^TZ}}], \mathbb{E}[\e^{4\bs{\beta_0^TZ}}], ||\alpha_0||_{\infty,\tau},\tau, f_0)\dfrac{\log(pn^k)}{n\sqrt{n}},
\end{align*}
where $s$ is the sparsiy index of $\bs{\beta_0}$ and 
\[
C(s,B,R,\E[\e^{\bs{\beta_0^TZ}}], ||\alpha_0||_{\infty,\tau},\tau)\quad  \text{and} \quad  \tilde{C}(s,B,|\bs{\beta_0}|_1,R,\mathbb{E}[\e^{\bs{\beta_0^TZ}}], \mathbb{E}[\e^{4\bs{\beta_0^TZ}}], ||\alpha_0||_{\infty,\tau},\tau, f_0)
\]
 are constants depending on the elements in brackets. 
Combining the previous bounds with Proposition \ref{propcomp}, we conclude that  Theorem \ref{th:IO2} is proved since 
\[
\mathbb{E}[||\hat\alpha^{\bs{\hat\beta}}_{\hat m^{\bs{\hat\beta}}}-\alpha^{\bs{\beta_0}}_m||^2_{det}]\leq \kappa_0\underset{m\in\mathcal{M}_n}{\inf}\{||\alpha_0-\alpha^{\bs{\beta_0}}_m||^2_{det}+2\pen(m)\}+\dfrac{C_1}{n}+C_2\dfrac{\log(pn)}{n},
\]
where $C_1$ and $C_2$ are constants depending on the sparsity index $s$ of $\bs{\beta_0}$, $B$, $|\bs{\beta_0}|_1$, $\mathbb{E}[\e^{\bs{\beta_0^TZ}}]$, $\mathbb{E}[\e^{4\bs{\beta_0^TZ}}]$,$||\alpha_0||_{\infty,\tau}$, $\tau$, ${f}_0$.

\qed
%%% SECTION 6 --- PROOFS OF THE TECHNICAL LEMMAS ---

\subsubsection{Proof of Corollary \ref{cor:Minimax}}

From Proposition \ref{prop:connection2} and the proof of Corollary 1 in \citet{CGG}, we deduce that 
\begin{align*}
\mathbb{E}[||\hat\alpha^{\bs{\hat\beta}}_{\hat m^{\bs{\hat\beta}}}-\alpha_0||^2_{2}]&\leq \dfrac{\e^{B|\bs{\beta_0}|_1}}{f_0}\mathbb{E}[||\hat\alpha^{\bs{\hat\beta}}_{\hat m^{\bs{\hat\beta}}}-\alpha_0||^2_{det}] \leq\tilde C_1\underset{m\in\mathcal{M}_n}{\inf}\Bigg\{D_m^{-2\gamma}+\dfrac{D_m}{n}\Bigg\}+\tilde C_2(s)\dfrac{\log(np)}{n},
\end{align*}
and since
\[
\underset{m\in\mathcal{M}_n}{\inf}\Bigg\{D_m^{-2\gamma}+\dfrac{D_m}{n}\Bigg\}=n^{-\frac{2\gamma}{2\gamma+1}},
\]
we finally get the corollary.
\qed

\subsection{Proofs of the technical propositions and lemmas}
\label{proof:technic2}

\subsubsection{Proof of Lemma \ref{lem:incl2}}
\label{proof:incl}

%Proofs of Lemma \ref{lem:incl1} and Lemma \ref{lem:incl2} rely on the same principle. We only detail the one that establishes Lemma \ref{lem:incl2} and we refer to \citet{CGG} for a similar proof of Lemma \ref{lem:incl1}.

%\paragraph{Proof of Lemma \ref{lem:incl1}. }
%Let $m\in\mathcal{M}_n$ be fixed and let $v$ be an eigenvalue of $\bs{G^{\bs{\beta_0}}_m}$. There exists $A_m\neq 0$ with coefficients $(a_j)_j$ such that $\bs{G^{\bs{\beta_0}}_m}A_m=vA_m$ and thus $A^T_m\bs{G^{\bs{\beta_0}}_m}A_m=vA^T_mA_m$. Now, take $h:=\sum_ja_j\varphi_j\in S_m$. We have $||h||^2_{rand}=A^T_m\bs{G^{\bs{\beta_0}}_m}A^T_m$ and $||h||^2_2=A^T_mA_m$.
%Thus, on $\Delta_1$ defined in (\ref{set:evenements2}) and from Proposition \ref{prop:connection2}:
%\[
%A^T_m\bs{G^{\bs{\beta_0}}_m}A^T_m=||h||^2_{rand}\geq \dfrac{1}{2}||h||^2_{det}\geq \dfrac{1}{2}f_0\e^{-B|\bs{\beta_0}|_1}||h||^2_2.
%\]
%Therefore, on $\Delta_1$, for all $m\in\mathcal{M}_n$, we have $\min\Sp(\bs{G^{\bs{\beta_0}}_m})\geq f_0\e^{-B|\bs{\beta_0}|_1}/2$. Moreover, on $\Omega$, we have $f_0\geq 2\hat{f}_0/3$ and $\max(\hat{f}_0/6,n^{-1/2})=\hat{f}_0/6$ for $n\geq 16/(f_0\e^{-B|\bs{\beta_0}|_1})^2$. \qed

%\paragraph{Proof of Lemma \ref{lem:incl2}. }
Let $m\in\mathcal{M}_n$ be fixed and let $v$ be an eigenvalue of $\bs{G^{\bs{\hat\beta}}_m}$. There exists $\bs{A_m}\neq 0$ with coefficients $(a_j)_j$ such that $\bs{G^{\bs{\hat\beta}}_m}\bs{A_m}=v\bs{A_m}$ and thus $\bs{A^T_m}\bs{G^{\bs{\hat\beta}}_m}\bs{A_m}=v\bs{A^T_m}\bs{A_m}$. Now, take $h:=\sum_ja_j\varphi_j\in S_m$. We have $||h||^2_{rand(\bs{\hat\beta})}=\bs{A^T_m}\bs{G^{\bs{\hat\beta}}_m}\bs{A^T_m}$ and $||h||^2_2=\bs{A^T_mA_m}$.
Thus, on $\Delta_1\cap\Delta_2$ defined in (\ref{set:evenements2}) and (\ref{set:delta2}) and from Proposition \ref{prop:connection2}:
\[
\bs{A^T_m\bs{G^{\bs{\hat\beta}}_m}A^T_m}=||h||^2_{rand(\bs{\hat\beta})}\geq \dfrac{1}{2}||h||^2_{rand}\geq \dfrac{1}{4}||h||^2_{det}\geq \dfrac{1}{4}f_0\e^{-B|\bs{\beta_0}|_1}||h||^2_2.
\]
Therefore, on $\Delta_1\cap\Delta_2$, for all $m\in\mathcal{M}_n$, we have $\min\Sp(\bs{G^{\bs{\hat\beta}}_m})\geq f_0\e^{-3BR}/4$. Moreover, on $\Omega$, we have $f_0\geq 2\hat{f}_0/3$ and $\max(\hat{f}_0\e^{-3BR}/6,n^{-1/2})=\hat{f}_0\e^{-3BR}/6$ for $n\geq 36/(\hat {f}_0\e^{-3BR})^2$, which is equivalent on $\Omega$ to choose $n\geq 16/(f_0\e^{-3BR})^2$. \qed

\subsubsection{Proof of Proposition \ref{propcomp}}
\label{proof:prop}

We have the following decomposition :
\begin{align*}
\mathbb{E}[||\hat\alpha^{\bs{\hat\beta}}_{\hat m^{\bs{\hat\beta}}}-\alpha^{\bs{\beta_0}}_m||^2_{det}\mathds{1}_{\aleph_k^c}]\leq &\mathbb{E}[||\hat\alpha^{\bs{\hat\beta}}_{\hat m^{\bs{\hat\beta}}}-\alpha^{\bs{\beta_0}}_m||^2_{det}\mathds{1}_{\Delta^c_1}]+\mathbb{E}[||\hat\alpha^{\bs{\hat\beta}}_{\hat m^{\bs{\hat\beta}}}-\alpha^{\bs{\beta_0}}_m||^2_{det}\mathds{1}_{\Delta^c_2}]\\
+&\mathbb{E}[||\hat\alpha^{\bs{\hat\beta}}_{\hat m^{\bs{\hat\beta}}}-\alpha^{\bs{\beta_0}}_m||^2_{det}\mathds{1}_{\Omega^c}]+\mathbb{E}[||\hat\alpha^{\bs{\hat\beta}}_{\hat m^{\bs{\hat\beta}}}-\alpha^{\bs{\beta_0}}_m||^2_{det}\mathds{1}_{(\Omega^k_H)^c}].
\end{align*}
We deduce that
\begin{align*}
\mathbb{E}[||\hat\alpha^{\bs{\hat\beta}}_{\hat m^{\bs{\hat\beta}}}-\alpha^{\bs{\beta_0}}_m||^2_{det}\mathds{1}_{\aleph_k^c}]\leq2 \Big(&\mathbb{E}[||\hat\alpha^{\bs{\hat\beta}}_{\hat m^{\bs{\hat\beta}}}-\alpha_0||^2_{det}\mathds{1}_{\Delta^c_1}]+\mathbb{E}[||\alpha^{\bs{\beta_0}}_m-\alpha_0||^2_{det}\mathds{1}_{\Delta^c_1}]\\
+&\mathbb{E}[||\hat\alpha^{\bs{\hat\beta}}_{\hat m^{\bs{\hat\beta}}}-\alpha_0||^2_{det}\mathds{1}_{\Delta^c_2}]+\mathbb{E}[||\alpha^{\bs{\beta_0}}_m-\alpha_0||^2_{det}\mathds{1}_{\Delta^c_2}]\\
+&\mathbb{E}[||\hat\alpha^{\bs{\hat\beta}}_{\hat m^{\bs{\hat\beta}}}-\alpha_0||^2_{det}\mathds{1}_{\Omega^c}]+\mathbb{E}[||\alpha^{\bs{\beta_0}}_m-\alpha_0||^2_{det}\mathds{1}_{\Omega^c}]\\
+&\mathbb{E}[||\hat\alpha^{\bs{\hat\beta}}_{\hat m^{\bs{\hat\beta}}}-\alpha_0||^2_{det}\mathds{1}_{(\Omega^k_H)^c}]+\mathbb{E}[||\alpha^{\bs{\beta_0}}_m-\alpha_0||^2_{det}\mathds{1}_{(\Omega^k_H)^c}]\Big).
\end{align*}
From definition (\ref{def:alpham}) of $\alpha^{\bs{\beta_0}}_m$ and Proposition \ref{prop:connection2}, we have $||\alpha^{\bs{\beta_0}}_m-\alpha_0||^2_{det}\leq||\alpha_0||^2_{det}\leq \mathbb{E}[\e^{\bs{\beta_0^TZ}}]||\alpha_0||^2_2$. From this relation and using Cauchy-Schwarz Inequality, we have 
\begin{align*}
\mathbb{E}[||\hat\alpha^{\bs{\hat\beta}}_{\hat m^{\bs{\hat\beta}}}-\alpha^{\bs{\beta_0}}_m||^2_{det}\mathds{1}_{\aleph_k^c}]\leq 4\mathbb{E}[\e^{\bs{\beta_0^TZ}}]\Big[\mathbb{E}^{1/2}(||\hat\alpha^{\bs{\hat\beta}}_{\hat m^{\bs{\hat\beta}}}||^4_2)\Big(\mathbb{P}^{1/2}(\Delta^c_1)+&\mathbb{P}^{1/2}(\Delta^c_2)\\
+\mathbb{P}^{1/2}(\Omega^c)+\mathbb{P}^{1/2}((\Omega^k_H)^c)\Big)+||\alpha_0||^2_2(\mathbb{P}(\Delta^c_1)+\mathbb{P}(\Delta^c_2)+&\mathbb{P}(\Omega^c)+\mathbb{P}((\Omega^k_H)^c))\Big].
\end{align*}
%From Assumption \ref{ass:assnot2}.\ref{ass:hatf0}, with $k=6$, we have immediately 
%\begin{equation}
%\label{Omega21}
%\mathbb{P}(\Omega^c)\leq C_0/n^6, \quad \forall n\geq n_0. 
%\end{equation}
From Assumption \ref{ass:hatf0}, Proposition \ref{prop:predbeta2}, Lemmas \ref{hatalpha}, \ref{Delta1} and \ref{Delta2} with $k=6$, we conclude that 
\begin{align*}
\mathbb{E}[||\hat\alpha^{\bs{\hat\beta}}_{\hat m^{\bs{\hat\beta}}}-\alpha^{\bs{\beta_0}}_m||^2_{det}\mathds{1}_{\aleph_k^c}]\leq& \hspace{0.1cm}2\mathbb{E}[\e^{\bs{\beta_0^TZ}}]\Bigg[\sqrt{C_bn^4}\Bigg(\sqrt{\dfrac{C^{(\Delta_1)}_6}{n^6}}+\sqrt{\dfrac{C^{(\Delta_2)}_6}{n^6}}+\sqrt{\dfrac{C_0}{n^6}}+\sqrt{\dfrac{c}{n^6}}\Bigg)\\
&+||\alpha_0||^2_2\Bigg(\dfrac{C^{(\Delta_1)}_6}{n^6}+\dfrac{C^{(\Delta_2)}_6}{n^6}+\dfrac{C_0}{n^6}+\dfrac{c}{n^6}\Bigg)\Bigg]\\
\leq &\hspace{0.1cm}\dfrac{\tilde c_1}{n},
\end{align*}
which ends the proof of Proposition \ref{propcomp}.
\qed

\subsubsection{Proof of Proposition \ref{espeta}}
\label{proof:espeta}
The proof is inspired from the paper of \cite{BCL10}.
If we denote $(\varphi_j)_{j\in\mathcal{K}_n}$ the orhonormal basis of the global nesting space $\mathcal{S}_n$ (see Assumption \ref{ass:model2}.\ref{ass:emboite}), since $\alpha$ belongs to $\mathcal{B}^{det}_n(0,1)\subset\mathcal{S}_n$, we can write $\alpha(t)=\sum_{j\in\mathcal{K}_n}a_j\varphi_j(t)$, with $\dim \mathcal{S}_n=\mathcal{D}_n=|\mathcal{K}_n|$. With this definition, we obtain
\begin{align*}
\eta_n(\alpha,\alpha^{\bs{\beta_0}}_m)=\displaystyle\sum_{j,j'}&a_ja_{j'}\dfrac{1}{n}\sumin\Big(\inttau\varphi_j(t)\alpha^{\bs{\beta_0}}_m(t)\expbeta Y_i(t)\diff t\inttau\varphi_{j'}\alpha^{\bs{\beta_0}}_m(t)\expbeta  Y_i(t)\diff t\\
&-\mathbb{E}\Big[\inttau\varphi_j(t)\alpha^{\bs{\beta_0}}_m(t)\expbeta Y_i(t)\diff t\inttau\varphi_{j'}\alpha^{\bs{\beta_0}}_m(t)\expbeta Y_i(t)\diff t\Big]\Big)\Big.
\end{align*}
For sake of simplicity, we introduce the notation
\[
A^{i}_{j,j'}= \inttau\varphi_j(t)\alpha^{\bs{\beta_0}}_m(t)\expbeta Y_i(t)\diff t\inttau\varphi_{j'}(t)\alpha^{\bs{\beta_0}}_m(t)\expbeta Y_i(t)\diff t.
\]
Applying the Cauchy-Schwarz Inequality, we get
\begin{align*}
|\eta_n(\alpha,\alpha^{\bs{\beta_0}}_m)|\leq \sqrt{\displaystyle\sum_{j,j'}a^2_ja^2_{j'}}\sqrt{\displaystyle\sum_{j,j'}\Big(\dfrac{1}{n}\sumin(A_{j,j'}^i-\mathbb{E}[A_{j,j'}^i])\Big)^2}.
\end{align*}
%With this notation, we have
%\[
%\eta_n(\alpha,\alpha^{\bs{\beta_0}}_m)^2\leq \dfrac{1}{n^2}\displaystyle\sum_{i,i'}\sum_{j,j'}a^2_ja^2_{j'}\sqrt{\sum_{j,j'}(A^i_{j,j'}-\mathbb{E}(A^i_{j,j'}))^2}\sqrt{\sum_{j,j'}(A^{i'}_{j,j'}-\mathbb{E}(A^{i'}_{j,j'}))^2}
%\]
From Proposition \ref{prop:connection2}, we have
\begin{align*}
\underset{\alpha\in\mathcal{B}^{det}_n(0,1)}{\sup}\eta_n(\alpha,\alpha^{\bs{\beta_0}}_m)^2&\leq\underset{(a_j),\sum_{j}a^2_j\leq 1}{\sup}\dfrac{1}{(\e^{-B|\bs{\beta_0}|_1}f_0)^2}\sum_{j,j'}a^2_ja^2_{j'}\sum_{j,j'}\Big(\dfrac{1}{n}\sumin (A^i_{j,j'}-\mathbb{E}[A^i_{j,j'}])\Big)^2\\
&\leq \dfrac{1}{(\e^{-B|\bs{\beta_0}|_1}f_0)^2}\sum_{j,j'}\Big(\dfrac{1}{n}\sumin(A^i_{j,j'}-\mathbb{E}[A^i_{j,j'}])\Big)^2.
\end{align*}
%For $i\neq i'$, $\sqrt{\sum_{j,j'}(A^i_{j,j'}-\mathbb{E}(A^i_{j,j'}))^2}$ and $\sqrt{\sum_{j,j'}(A^{i'}_{j,j'}-\mathbb{E}(A^{i'}_{j,j'}))^2}$ are independant, hence, taking the expectation, we obtain that $\mathbb{E}[{\sup}_{\alpha\in\mathcal{B}^{det}_n(0,1)}\eta_n(\alpha,\alpha^{\bs{\beta_0}}_m)^2]$ is less than
%\begin{align*}
%\dfrac{1}{(\e^{-B|\bs{\beta_0}|}f_0)^2}\Bigg(\dfrac{1}{n^2}\sum_{i\neq i'}&\mathbb{E}\Bigg[\sqrt{\sum_{j,j'}(A^i_{j,j'}-\mathbb{E}(A^i_{j,j'}))^2}\Bigg]\mathbb{E}\Bigg[\sqrt{\sum_{j,j'}(A^{i'}_{j,j'}-\mathbb{E}(A^{i'}_{j,j'}))^2}\Bigg]\\
%+\sum_{i=1}^{n}&\mathbb{E}\Bigg[\sum_{j,j'}(A^i_{j,j'}-\mathbb{E}(A^i_{j,j'}))^2\Bigg]\Bigg).
%\end{align*}
%It follows that
%\begin{align*}
Taking the expectation, it follows that 
\begin{align*}
\mathbb{E}\Bigg[\underset{\alpha\in\mathcal{B}^{det}_n(0,1)}{\sup}\eta_n(\alpha,\alpha^{\bs{\beta_0}}_m)^2\Bigg]&\leq  \dfrac{1}{(\e^{-B|\bs{\beta_0}|_1}f_0)^2}\sum_{j,j'}\var\Bigg[\dfrac{1}{n}\sumin A^i_{j,j'}\Bigg]\\
&\leq  \dfrac{1}{(\e^{-B|\bs{\beta_0}|_1}f_0)^2}\sum_{j,j'}\dfrac{1}{n}\E\Big[(A^1_{j,j'})^2\Big].
\end{align*}
%and with an unrefined bound, we obtain
%\begin{align*}
%\leq \dfrac{1}{(\e^{-B|\bs{\beta_0}|}f_0)^2}\Bigg(\dfrac{1}{n^2}&\sum_{i=1}^{n} \Bigg(\mathbb{E}\Bigg[{\sum_{j,j'}(A^i_{j,j'}-\mathbb{E}(A^i_{j,j'}))}\Bigg]\Bigg)^2\\
%+&\sum_{i=1}^{n}\mathbb{E}\Bigg[\sum_{j,j'}(A^i_{j,j'}-\mathbb{E}(A^i_{j,j'}))^2\Bigg]\Bigg).
%\end{align*}
%Since $\Big(\mathbb{E}\Big[{\sum_{j,j'}\Big(A^i_{j,j'}-\mathbb{E}(A^i_{j,j'})\Big)}\Big]\Big)^2=0$, we obtain that $\mathbb{E}[{\sup}_{\alpha\in\mathcal{B}^{det}_n(0,1)}\eta_n(\alpha,\alpha^{\bs{\beta_0}}_m)^2]$ is less than
%\begin{align*}
%&\dfrac{1}{(\e^{-B|\bs{\beta_0}|}f_0)^2}\Bigg(\dfrac{1}{n^2}\sum_{i=1}^{n} \sum_{j,j'}\mathbb{E}\Bigg[(A^i_{j,j'})^2-2A^i_{j,j'}\mathbb{E}(A^i_{j,j'})+(\mathbb{E}(A^i_{j,j'}))^2\Bigg]\Bigg)\\
%\leq &\dfrac{1}{(\e^{-B|\bs{\beta_0}|}f_0)^2}\dfrac{1}{n^2}\sum_{i=1}^{n} \sum_{j,j'}\mathbb{E}\Big[(A^i_{j,j'})^2\Big]-\Big(\mathbb{E}\Big[A^i_{j,j'}\Big]\Big)^2\\
%\leq &\dfrac{1}{(\e^{-B|\bs{\beta_0}|}f_0)^2}\dfrac{1}{n} \sum_{j,j'}\mathbb{E}\Bigg[\Bigg(\inttau\varphi_j(t)\alpha^{\bs{\beta_0}}_m(t)\e^{\bs{\beta_0^TZ}}Y(t)\diff t\Bigg)^2\Bigg(\inttau\varphi_{j'}(t)\alpha^{\bs{\beta_0}}_m(t)\e^{\bs{\beta_0^TZ}}Y(t)\diff t\Bigg)^2\Bigg].
%\end{align*}
Thus, from the definition of $A^1_{j,j'}$, we obtain that $\mathbb{E}[{\sup}_{\alpha\in\mathcal{B}^{det}_n(0,1)}\eta_n(\alpha,\alpha^{\bs{\beta_0}}_m)^2]$ is less than
\begin{align*}
\dfrac{1}{(\e^{-B|\bs{\beta_0}|}f_0)^2}\dfrac{1}{n} \sum_{j,j'}\mathbb{E}\Bigg[\Bigg(\inttau\varphi_j(t)\alpha^{\bs{\beta_0}}_m(t)\e^{\bs{\beta_0^TZ}}Y(t)\diff t\Bigg)^2\Bigg(\inttau\varphi_{j'}(t)\alpha^{\bs{\beta_0}}_m(t)\e^{\bs{\beta_0^TZ}}Y(t)\diff t\Bigg)^2\Bigg].
\end{align*}
From \citet{BCL10} p.301, Equation (2.7), we have 
\[
\displaystyle\sum_{j\in \mathcal{K}_n}\left(\inttau\varphi_j(t)\alpha^{\bs{\beta_0}}_m(t)\e^{\bs{\beta_0^TZ}}Y(t)\diff t\right)^2\leq \inttau (\alpha^{\bs{\beta_0}}_m(t)\e^{\bs{\beta_0^TZ}}Y(t))^2\diff t\leq \e^{2\bs{\beta_0^TZ}}||\alpha^{\bs{\beta_0}}_m||^2_2.
\]
From this inequality, we obtain
\begin{align*}
\mathbb{E}\left[\underset{\alpha\in\mathcal{B}^{det}_n(0,1)}{\sup}\eta_n(\alpha,\alpha^{\bs{\beta_0}}_m)^2\right]\leq \dfrac{ \mathbb{E}[\e^{4\bs{\beta_0^TZ}}]||\alpha^{\bs{\beta_0}}_m||^4_2}{(\e^{-B|\bs{\beta_0}|_1}f_0)^2}\dfrac{1}{n}.\qed
\end{align*}

\subsubsection{Proof of Lemma \ref{hatalpha}}
\label{proof:hatalpha}
%Let $\bs{\beta}\in\mathbb{R}^p$, so that $\bs{\beta}$ can be either $\bs{\beta_0}$ or $\bs{\hat\beta}$. 
%For sake of simplicity in the notations, we omit the parameter $\bs{\beta}$ in $\hat m^{\bs{\beta}}$. 
From Assumption \ref{ass:betaball2}, we recall that $|\bs{\hat\beta}-\bs{\beta_0}|_1\leq 2R$.
On $\mathcal{\hat{H}}^{\bs{\hat\beta}}_{\hat m^{\bs{\hat\beta}}}$, we have
%Under Assumption \ref{eigenvalue}, we have
\begin{align*}
||\hat\alpha^{\bs{\hat\beta}}_{\hat m^{\bs{\hat\beta}}}||^2_2&=\displaystyle\sum_{j\in J_{\hat m^{\bs{\hat\beta}}}}(\hat{a}^{\hat m^{\bs{\hat\beta}}}_j)^2=||\bs{A_{\hat m^{\bs{\hat\beta}}}}||^2_2=||(\bs{G^{\bs{\hat\beta}}_{\hat m^{\bs{\hat\beta}}}})^{-1}\bs{\Gamma_{\hat m^{\bs{\hat\beta}}}}||^2_2\\
&\leq (\min\Sp(\bs{G^{\hat\beta}_{\hat m^{\bs{\hat\beta}}}}))^{-2}||\bs{\Gamma_{\hat m^{\bs{\hat\beta}}}}||^2_2\\
&\leq\min\left(\dfrac{36}{\hat{f}^2_0 \e^{-2B|\bs\beta_0|_1-2B|\bs{\beta_0}-\bs{\hat\beta}|_1}},n\right)\displaystyle\sum_{j\in J_{\hat{m}^{\bs{\hat\beta}}}}\left(\dfrac{1}{n}\sumin\inttau\varphi_j(t)\diff N_i(t)\right)^2\\
&\leq  \min\left(\dfrac{36}{\hat{f}^2_0 \e^{-2B|\bs\beta_0|_1-4BR}},n\right)\dfrac{1}{n}\sumin\displaystyle\sum_{j\in J_{\hat{m}^{\bs{\hat\beta}}}}\left(\inttau\varphi_j(t)\diff N_i(t)\right)^2.
\end{align*} 
So we have
\begin{align*}
||\hat\alpha^{\bs{\hat\beta}}_{\hat m^{\bs{\hat\beta}}}||^4_2\leq n^2\dfrac{1}{n}\sumin\left(\displaystyle\sum_{j\in J_{\hat{m}^{\bs{\hat\beta}}}}\left(\inttau\varphi_j(t)\diff N_i(t)\right)^2\right)^2\leq n^2\dfrac{1}{n}\sumin\left(\displaystyle\sum_{j\in \mathcal{K}_n}\left(\inttau\varphi_j(t)\diff N_i(t)\right)^2\right)^2,
\end{align*}
where $\mathcal{K}_n$ is a set of indices of the global nesting space $\mathcal{S}_n$, defined in Assumption \ref{ass:model2}.\ref{ass:emboite}, and $\dim \mathcal{S}_n=\mathcal{D}_n=|\mathcal{K}_n|$.
Thus, we deduce that
\begin{align*}
||\hat\alpha^{\bs{\hat\beta}}_{\hat m^{\bs{\hat\beta}}}||^4_2\leq n^2 \mathcal{D}_{n}\dfrac{1}{n}\sumin\displaystyle\sum_{j\in \mathcal{K}_n}\left(\inttau\varphi_j(t)\diff N_i(t)\right)^4.
\end{align*}
Now, 
\begin{align*}
\mathbb{E}\left[\dfrac{1}{n}\sumin\displaystyle\sum_{j\in  \mathcal{K}_n}\left(\inttau\varphi_j(t)\diff N_i(t)\right)^4\right]&\leq \dfrac{2^3}{n}\sumin\displaystyle\sum_{j\in  \mathcal{K}_n}\mathbb{E}\left[\left(\inttau\varphi_j(t)\diff M_i(t)\right)^4\right]\\
&+\dfrac{2^3}{n}\sumin\displaystyle\sum_{j\in  \mathcal{K}_n}\mathbb{E}\left[\left(\inttau\varphi_j(t)\alpha_0(t)\expbeta Y_i(t)\diff t\right)^4\right].
\end{align*}
Using the B\"urkholder Inequality (see \citet{Liptser89}), we get 
\begin{align*}
\mathbb{E}\left[\dfrac{1}{n}\sumin\displaystyle\sum_{j\in  \mathcal{K}_n}\left(\inttau\varphi_j(t)\diff M_i(t)\right)^4\right]&\leq \kappa_b\dfrac{1}{n}\sumin\displaystyle\sum_{j\in  \mathcal{K}_n}\mathbb{E}\left[\left(\inttau\varphi^2_j(t)\diff N_i(t)\right)^2\right]\\
&\leq \kappa_b\dfrac{1}{n}\sumin\displaystyle\sum_{j\in  \mathcal{K}_n}\mathbb{E}\left[N_i(\tau)\displaystyle\sum_{s:\Delta N_i\neq 0}\varphi^4_j(s)\right]\\
&\leq \kappa_b\dfrac{1}{n}\sumin\mathbb{E}\left[N_i(\tau)\displaystyle\sum_{s:\Delta N_i\neq 0}\displaystyle\sum_{j\in  \mathcal{K}_n}\varphi^4_j(s)\right],\\
\end{align*}
which is finally bounded from Assumption \ref{ass:model2}.\ref{ass:norminf2} by
\begin{align*}
\mathbb{E}\left[\dfrac{1}{n}\sumin\displaystyle\sum_{j\in  \mathcal{K}_n}\left(\inttau\varphi_j(t)\diff M_i(t)\right)^4\right]&\leq \kappa_b\phi^2\mathcal{D}_n^2\dfrac{1}{n}\sumin\mathbb{E}\left[N_i(\tau)\displaystyle\sum_{s:\Delta N_i\neq 0}1\right]\\
&\leq \kappa_b\phi^2\mathcal{D}_n^2\mathbb{E}[N_1(\tau)^2].
\end{align*}
Then, we can write that 
\begin{align*}
[N_1(\tau)]^2&=\left[M_1(\tau)+\inttau\alpha_0(t)\e^{\bs{\beta_0^TZ}} Y(t)\diff t\right]^2\\
&\leq 2(M_1(\tau))^2+2\left(\inttau\alpha_0(t)\e^{\bs{\beta_0^TZ}} Y(t)\diff t\right)^2, 
\end{align*}
and 
\[
\mathbb{E}[(M_1(\tau))^2]\leq \mathbb{E}\left[\inttau\alpha_0(t)\e^{\bs{\beta_0^TZ}} Y(t)\diff t\right]\leq \tau||\alpha_0||_{\infty,\tau}\mathbb{E}[\e^{\bs{\beta_0^TZ}}],
\]
so that 
\[
\mathbb{E}[(N_1(\tau))^2]\leq 2||\alpha_0||_{\infty,\tau}\tau\mathbb{E}[\e^{\bs{\beta_0^TZ}}]+2||\alpha_0||^2_{\infty,\tau}(\mathbb{E}[\e^{\bs{\beta_0^TZ}}])^2\tau^2.
\]
So, by using Cauchy-Schwarz Inequality, we obtain
\begin{multline*}
\mathbb{E}\Bigg[\dfrac{1}{n}\sumin\displaystyle\sum_{j\in \mathcal{K}_n}\Bigg(\inttau\phi_j(t)\diff N_i(t)\Bigg)^4\Bigg]\\
\leq 8\kappa_b\phi^2\mathcal{D}_n^2\mathbb{E}[(N_1(\tau))^2]+8\displaystyle\sum_{j\in \mathcal{K}_n}\mathbb{E}\Bigg[\Bigg(\inttau\varphi_j(t)\alpha_0(t)\e^{\bs{\beta_0^TZ}} Y(t)\diff t\Bigg)^4\Bigg]\\
%\leq &8\kappa_b\phi^2D_{\hat m}^2\mathbb{E}[(N_1(\tau))^2]+8||\alpha_0||^4_{\infty,\tau}(\mathbb{E}[\e^{2\bs{\beta_0^TZ}}])^2\tau^2\displaystyle\sum_{j\in J_{\hat m}}\Bigg(\inttau\varphi^2_j(t)\diff t\Bigg)^2\\
\leq 8\kappa_b\phi^2\mathcal{D}_n^2\mathbb{E}[(N_1(\tau))^2]+8||\alpha_0||^4_{\infty,\tau}\mathbb{E}[\e^{4\bs{\beta_0^TZ}}]\tau^2\mathcal{D}_n.
\end{multline*}
Eventually, under Assumption \ref{ass:model2}.\ref{ass:Dm2}, we get 
\begin{align*}
\mathbb{E}[||\hat\alpha^{\bs{\hat\beta}}_{\hat m^{\bs{\hat\beta}}}||^4_2]&\leq n^2\mathcal{D}_n\Big[8\kappa_b\phi^2\mathcal{D}_n^2\Big(2||\alpha_0||_{\infty\tau}\tau\mathbb{E}[\e^{\bs{\beta_0^TZ}}]+2||\alpha_0||^2_{\infty,\tau}(\mathbb{E}[\e^{\bs{\beta_0^TZ}}])^2\tau^2\Big)\\
&\hspace{1.2cm}+8||\alpha_0||^4_{\infty,\tau}\mathbb{E}[\e^{4\bs{\beta_0^TZ}}]\tau^2\mathcal{D}_n\Big]\\
&\leq C_b n^2\mathcal{D}_n^3\\
&\leq C_b n^4,
\end{align*}
where $C_b$ is a constant that depends on $\kappa_b$, $||\alpha_0||_{\infty,\tau}$, $\tau$, $\mathbb{E}[\e^{\bs{\beta_0^TZ}}]$ and $\mathbb{E}[\e^{4\bs{\beta_0^TZ}}]$ and on the choice of the basis.
\qed

\subsubsection{Proof of Lemma \ref{Delta1} }
\label{proof:Delta1}
The event $\Delta_1$ defined by (\ref{set:evenements2}) can be rewritten as
\[
\Delta_1=\left\{\omega\in\Omega, \forall \alpha\in\mathcal{S}_n\backslash\{0\} : \left|\dfrac{||\alpha||^2_{rand(\omega)}}{||\alpha||^2_{det}}-1\right|\leq \dfrac{1}{2}\right\},
\]
and consider
\begin{equation}
\label{vartheta}
\vartheta_n(\alpha)=\dfrac{1}{n}\sumin\inttau\Big(\alpha(t)\expbeta Y_i(t)-\mathbb{E}[\alpha(t)\expbeta Y_i(t)]\Big)\diff t=||\sqrt{\alpha}||^2_{rand}-||\sqrt{\alpha}||^2_{det}.
\end{equation}
If $\omega\in(\Delta_1)^c$, then there exists $\alpha$ (which can depend on $\omega$) such that 
\[
\Bigg|\dfrac{||\alpha||^2_{rand(\omega)}}{||\alpha||^2_{det}}-1\Bigg|>\dfrac{1}{2}.
\]
Taking $\gamma=\alpha/||\alpha||^2_{det}$, we have that 
\[
 \gamma\in\mathcal{S}_n\backslash\{0\}, \quad ||\gamma||^2_{det}=1, \quad \text{ and } \quad
 |||\gamma||^2_{rand(\omega)}-1|>\dfrac{1}{2}.
 \]
So, if $\omega\in(\Delta_1)^c$, then 
\[
\omega\in\Bigg\{\omega\in\Omega : \underset{\gamma\in\mathcal{S}_n\backslash\{0\}, ||\gamma||^2_{det}=1}{\sup}|||\gamma||^2_{rand(\omega)}-1|>\dfrac{1}{2}\Bigg\}
\]
From this, we deduce that,
\[
\mathbb{P}((\Delta_1)^c)\leq \mathbb{P}\left(\underset{\alpha\in\mathcal{B}^{det}_{n}(0,1)}{\sup}|\vartheta_n(\alpha^2)|>1-\dfrac{1}{\rho_1}\right),
\]
where $\mathcal{B}^{det}_{n}(0,1)$ is defined by (\ref{eq:Bn01}).
%=\{\alpha\in\mathcal{S}_n : ||\alpha||_{det}\leq 1\}$.
Since $\alpha\in\mathcal{B}^{det}_{n}(0,1)\subset\mathcal{S}_n$, then we can write $\alpha(t)=\sum_{j \in \mathcal{K}_n}a^m_j\varphi_j(t)$, where $\mathcal{K}_n$ is a set of indices of $\mathcal{S}_n$ and $\dim\mathcal{S}_n=\mathcal{D}_n=|\mathcal{K}_n|$. With this notation, we have
\[
\vartheta_n(\alpha^2)=\displaystyle\sum_{j,k}a_ja_k\vartheta_n(\varphi_j\varphi_k).
\]
From Proposition \ref{prop:connection2}, we have
\[
\underset{\alpha\in\mathbb{B}^{det}_{n}(0,1)}{\sup}|\vartheta_n(\alpha^2)|\leq \dfrac{1}{f_0\e^{-B|\bs{\beta_0}|_1}}\underset{(a_j),\sum_{j\in \mathcal{K}_n}a^2_j\leq 1}{\sup}\Big|\displaystyle\sum_{j,k}a_ja_k\vartheta_n(\varphi_j\varphi_k)\Big|.
\]
Let consider the process $(U^{(j,k)}_i)$ defined by
\[
U^{(j,k)}_i=\inttau\varphi_j(t)\varphi_k(t)\expbeta Y_i(t)\diff t,
\]
We have $|U^{(j,k)}_i|\leq \e^{B|\bs\beta_0|_1}$ and from Cauchy-Schwarz Inequality, we have
\[
(U^{(j,k)}_i)^2\leq\e^{2B|\bs{\beta_0}|_1}\inttau\varphi^2_j(t)\diff t\inttau\varphi^2_k(t)\diff t\leq \e^{2B|\bs{\beta_0}|_1}.
\]
We can apply the standard Bernstein Inequality (see \cite{Massart}) to the process $(U^{(j,k)}_i)$, and we obtain 
\begin{equation}
\label{BernsVn}
\mathbb{P}\Big(|\vartheta_n(\varphi_j\varphi_k)|\geq \e^{B|\bs\beta_0|_1}x+\sqrt{2\e^{2B|\bs{\beta_0}|_1}x}\Big)\leq 2\e^{-nx}.
\end{equation}
Let introduce 
\[
\Theta:=\{\forall j,k, |\vartheta_n(\varphi_j\varphi_k)|\leq \e^{B|\bs{\beta_0}|_1}x+\e^{B|\bs{\beta_0}|_1}\sqrt{2x}\} \quad \mbox{and} \quad x:=\dfrac{f^2_0\e^{-2B|\bs\beta_0|_1}}{16\mathcal{D}_n^2\e^{2B|\bs{\beta_0}|_1}}.
\]
On $\Theta$, we can write that ${\sup}_{\alpha\in\mathcal{B}^{det}_{n}(0,1)}|\vartheta_n(\alpha^2)|$ is less than
\begin{align*}
&\dfrac{1}{f_0\e^{-B|\bs\beta_0|_1}}\underset{(a_j),\sum_{j\in \mathcal{K}_n}a^2_j\leq 1}{\sup}\displaystyle\sum_{j,k}|a_ja_k|(\e^{B|\bs\beta_0|_1}x+\e^{B|\bs\beta_0|_1}\sqrt{2x})\\
\leq &\dfrac{1}{f_0\e^{-B|\bs\beta_0|_1}}\underset{(a_j),\sum_{j\in \mathcal{K}_n}a^2_j\leq 1}{\sup}\Big(\displaystyle\sum_{j}|a_j|\Big)^2(\e^{B|\bs\beta_0|_1}x+\e^{B|\bs\beta_0|_1}\sqrt{2x}),\\
\end{align*}
which is less than
\begin{align}
\nonumber
\leq&\dfrac{1}{f_0\e^{-B|\bs\beta_0|_1}}D_m\left(\dfrac{\e^{B|\bs\beta_0|_1}f^2_0\e^{-2B|\bs\beta_0|_1}}{16\mathcal{D}_n^2\e^{2B|\bs\beta_0|_1}}+\dfrac{\e^{B|\bs\beta_0|_1}\sqrt{2}f_0\e^{-B|\bs\beta_0|_1}}{4\mathcal{D}_n\e^{B|\bs\beta_0|_1}}\right)\\
\nonumber
\leq&\dfrac{1}{2}\left(\dfrac{1}{8}\dfrac{f_0}{\e^{2B|\bs{\beta_0}|_1}\mathcal{D}_n}+\dfrac{1}{\sqrt{2}}\right)\\
\nonumber
\leq&\dfrac{1}{2}\left(\dfrac{1}{4}+\dfrac{1}{\sqrt{2}}\right)\\
\label{Theta}
\leq& \dfrac{1}{2}.
\end{align}
From Inequality (\ref{Theta}), we deduce that $\mathbb{P}((\Delta_1)^c)\leq \mathbb{P}(\Theta^c)$. So using Inequality (\ref{BernsVn}), we can conclude that
\begin{align*}
\mathbb{P}((\Delta_1)^c)&\leq \displaystyle\sum_{j,k}\mathbb{P}\Big(|\vartheta_n(\varphi_j\varphi_k)|>\e^{B|\bs\beta_0|_1}x+\e^{B|\bs\beta_0|_1}\sqrt{2x}\Big)\\
&\leq 2\mathcal{D}_n^2\exp\left(-\dfrac{nf^2_0\e^{-2B|\bs\beta_0|_1}}{16\mathcal{D}_n^2\e^{2B|\bs{\beta_0}|_1}}\right)\\
&\leq 2n\exp\left(-\dfrac{f^2_0}{16\e^{4B|\bs\beta_0|_1}}\dfrac{n}{\mathcal{D}_n^2}\right)\\
&\leq 2n\exp\left(-\dfrac{f^2_0}{16\e^{4B|\bs\beta_0|_1}}\log n\right)\\
&\leq \dfrac{C^{\Delta_1}_k}{n^k}, \quad \forall k\geq 1,
\end{align*}
as $\mathcal{D}_n\leq \sqrt{n}/\log n$ from Assumption \ref{ass:model2}.\ref{ass:emboite}, which ends the proof of Lemma \ref{Delta1} with $C^{\Delta_1}_k$ a constant depending on $\rho_1$, $f_0$, $B$ and $|\bs{\beta}_0|_1$.
\qed

\subsubsection{Proof of Lemma  \ref{Delta2}}
\label{proof:Delta2}
For $\rho_2\geq 1$, let define
\[
\Delta^{\rho_2}_2=\left\{\forall \alpha\in\mathcal{S}_n :\left|\dfrac{||\alpha||^2_{\rand}}{||\alpha||^2_{rand}}-1\right|\leq 1-\dfrac{1}{\rho_2}\right\}.
\]
Let consider 
\[
\tilde\vartheta_n(\alpha)=\dfrac{1}{n}\sumin\inttau(\alpha(t)\exphatbeta Y_i(t)-\alpha(t)\expbeta Y_i(t))\diff t=||\sqrt{\alpha}||^2_{\rand}-||\sqrt{\alpha}||^2_{rand}.
\]
Following the same approach as in the proof of Lemma \ref{Delta1}, we have
\begin{align}
\label{eq:delta2}
\mathbb{P}((\Delta^{\rho_2}_2)^c)\leq \mathbb{P}\left(\underset{\alpha\in\mathcal{B}^{det}_{n}(0,1)}{\sup}|\tilde\vartheta_n(\alpha^2)|>1-\dfrac{1}{\rho_2}\right),
\end{align}
where $\mathcal{B}^{det}_{n}(0,1)=\{\alpha\in\mathcal{S}_n : ||\alpha||_{det}\leq 1\}$.
The process $\tilde\vartheta_n(\alpha^2)$ is bounded by
\[
|\tilde\vartheta_n(\alpha^2)|\leq B\e^{B|\bs\beta_0|_1}\e^{2BR}|\bs{\hat\beta}-\bs{\beta_0}|_1||\alpha||^2_2\leq |\bs{\hat\beta}-\bs{\beta_0}|_1\dfrac{ B\e^{B|\bs\beta_0|_1}\e^{2BR}}{f_0\e^{-B|\bs{\beta_0}|_1}}||\alpha||^2_{det}.
\]
So we get 
\[
\underset{\alpha\in\mathcal{B}_{\mathcal{S}_n}^{det}(0,1)}{\sup}|\tilde\vartheta_n(\alpha^2)|\leq|\bs{\hat\beta}-\bs{\beta_0}|_1\dfrac{ B\e^{2B|\bs\beta_0|_1}\e^{2BR}}{f_0}. 
\]
From Proposition \ref{prop:predbeta2}, we have with probability larger than $1-cn^{-k}$ 
\[
|\bs{\hat\beta}-\bs{\beta_0}|_1\leq C(s)\sqrt{\dfrac{\log (pn^k)}{n}}.
\] 
Then we have with probability larger than $1-cn^{-k}$ 
\[
\underset{\alpha\in\mathcal{B}_{\mathcal{S}_n}^{det}(0,1)}{\sup}|\tilde\vartheta_n(\alpha^2)|\leq C(s)\sqrt{\dfrac{\log(pn^k)}{n}}\dfrac{B\e^{2B|\bs\beta_0|_1}\e^{2BR}}{f_0}.
\]
%We also have that 
%\[
%\mathbb{P}((\Delta^{\rho_2}_2)^c)\leq \mathbb{P}\Big(\underset{\alpha\in\mathcal{B}_{\mathcal{S}_n}^{det}(0,1)}{\sup}|\tilde\vartheta_n(\alpha^2)|>1-1/\rho_2\Big).
%\]
Thus, by taking $1-1/\rho_2=C(s)\sqrt{\dfrac{\log(pn^k)}{n}}\dfrac{B\e^{2B|\bs\beta_0|_1}\e^{2BR}}{f_0}$ in (\ref{eq:delta2}), we obtain
\[
\mathbb{P}((\Delta^{\rho_2}_2)^c)\leq cn^{-k}.
\]
From Assumption \ref{ass:Relationpn}, we deduce that for $n$ large enough, 
\[
1-\dfrac{1}{\rho_2}<\dfrac{1}{2},
\]
so that $\Delta_2$ defined by (\ref{set:delta2}) verifies $\mathbb{P}((\Delta_2)^c)\leq\mathbb{P}((\Delta^{\rho_2}_2)^c)\leq C^{(\Delta_2)}_kn^{-k}$, with $C^{(\Delta_2)}_k=c>0$.
%We also have $\mathbb{P}((\Delta^{\rho_2}_2)^c\cap(\Omega^k_H)^c)\leq \mathbb{P}((\Omega^k_H)^c)\leq cn^{-k}$ from Lemma \ref{lem:OmegaH} for a constant $c>0$ depending on $\tau$, $||\alpha_0||_{\infty,\tau}$ and $\mathbb{E}[\e^{\bs{\beta_0^T}}]$. 
%So finally, we obtain $\mathbb{P}((\Delta^{\rho_2}_2)^c)\leq C^{(\Delta_2)}_kn^{-k}$, with $C^{(\Delta_2)}_k=c>0$. 
%We have shown in Appendix \ref{ann:Huang} that the constant $c$ depends on $\tau$, $||\alpha_0||_{\infty,\tau}$ and $\mathbb{E}[\e^{\bs{\beta_0^TZ}}]$. 
\qed

\appendix

\section{Prediction result on the Lasso estimator ${\hat\beta}$ of ${\beta_0}$ for unbounded counting processes}
%\sectionmark{Prediction result on the Lasso estimator $\bs{\hat\beta}$ of $\bs{\beta_0}$ for $N_i$ unbounded}
\label{ann:Huang}

To obtain a non-asymptotic prediction bound on the Lasso estimator $\bs{\hat\beta}$ of the regression parameter in the Cox model, we rely on Theorem 3.1 of \citet{Huang2013}, that we recall here.

Let consider the classical Lasso estimator $\bs{\hat\beta}$ defined by (\ref{eq:betaL2}) when $p\gg n$.

We define $\bs{\dot l^*_n(\bs{\beta})}=(\dot l^*_{n,1}(\bs{\beta}),...,\dot l^*_{n,p}(\bs{\beta}))^T=\partial l^*_n(\bs{\beta})/\partial\bs{\beta}$ the gradient of the Cox partial log-likelihood $l^*_n(\bs{\beta})$ defined by (\ref{eq:PartialCoxLikelihood2}) and $\bs{\ddot l^*_n(\bs{\beta})}=\partial^2 l^*_n(\bs{\beta})/\partial\bs{\beta}\partial\bs{\beta}^T$ the Hessian matrix.

%Our results are deeply related to the result of \citet{Huang2013} obtained for the classical Lasso procedure (\ref{eq:betaL2}). 

Let us now describe the result of \citet{Huang2013}, on which we rely for our study, starting with the notations.
Let $\mathcal{O}=\{j : {\beta_0}_j\neq 0\}$, $\mathcal{O}^c=\{j : {\beta_0}_j=0\}$ and $s=|\mathcal{O}|$ the cardinality of $\mathcal{O}$. For any $\xi>1$, we define the cone    
\[
\mathcal{C}(\xi,\mathcal{O})=\{\bs{b}\in\mathbb{R}^p: |\bs{b}_{\mathcal{O}^c}|_1\leq \xi|\bs{b}_{\mathcal{O}}|_1\}.
\]
For this cone, let us define the following condition:
\[
0<\kappa(\xi,\mathcal{O})=\underset{0\neq\bs{b}\in\mathcal{C}(\xi,\mathcal{O})}{\inf}\dfrac{s^{1/2}(\bs{b\ddot l^*_n(\bs{\beta_0})b})^{1/2}}{|\bs{b}_{\mathcal{O}}|_1}.
\]
This term corresponds to the compatibility factor introduced by \citet{GeerRT}.  It is one of the classical condition used  to obtain non-asymptotic oracle inequalities. See also \citet{BvdG} for more details about this compatibility factor and the comparison of this criterion with other assumptions such as the Restricted Eigenvalue condition among other.
 
%We introduce the following assumption concerning the covariates :
%\begin{assumption}
%\label{ass:Zij2}
%There exists a positive constant $B$ such that for all $i\in\{1,...,n\}$ and $j\in\{1,...,p\}$, 
%\[|Z_{i,j}|\leq B.\]
%\end{assumption}
%This is a classical assumption (see Huang et al. \cite{Huang2013} and Bradic and Song \cite{bradic2012}).

With these notations, we can state the following theorem established by \citet{Huang2013}.

\begin{theorem}[\citet{Huang2013}]
\label{th:Huangann2}
Let $k>0$ and $\nu=B(\xi+1)s\Gamma_{n,k}/\{2\kappa^2(\xi,\mathcal{O})\}$. Suppose Assumption \ref{ass:assnot2}.\ref{ass:Zij2} holds and $\nu\leq 1/\e$. Then, on the event 
\begin{align}
\label{set:OmegaHann2}
\widetilde\Omega^k_H=\left\{|\bs{\dot l^*_n(\bs{\beta_0})}|_{\infty} \leq \dfrac{\xi-1}{\xi+1}\Gamma_{n,k}\right\}, \quad\mbox{ with } \quad \Gamma_{n,k}=C_0B\dfrac{\xi+1}{\xi-1}\sqrt{2\dfrac{\log(pn^k)}{n}},
\end{align}
we have
\[
|\bs{\hat\beta}-\bs{\beta_0}|_1\leq \dfrac{\e^{\eta}(\xi+1)s}{2\kappa^2(\xi,\mathcal{O})}\Gamma_{n,k},
\] 
where $\eta\leq 1$ is the smaller solution of $\eta\e^{-\eta}=\nu$ and $C_0>\sqrt{\tau||\alpha_0||_{\infty,\tau}\E[\e^{\bs{\beta_0^TZ}}]}$.
\end{theorem}
We refer to \citet{Huang2013} for the proof of Theorem \ref{th:Huangann2}. \citet{Huang2013} have calculated the probability of $\widetilde\Omega^k_H$ only in the case where $\max_{1\leq i\leq n}|N_i(\tau)|<+\infty$. We extend the result to the unbounded case in the following lemma.
%The following lemma gives the probability of the event $\widetilde\Omega^k_H$ defined by (\ref{set:OmegaHann2}). 
\begin{lemma}
\label{lem:OmegaHann}
Let consider, for $k>0$, the event $\widetilde\Omega^k_H$ defined by (\ref{set:OmegaHann2}).
%\begin{equation}
%\label{set:OmegaH2}
%\Omega^k_H=\left\{|\dot\ell(\bs{\beta_0})|_{\infty} \leq \dfrac{\xi-1}{\xi+1}\mu_k\right\}, \quad \mbox{with }  \mu_k=\dfrac{\xi+1}{\xi-1}B\sqrt{2\dfrac{\log(pn^k)}{n}}.
%\end{equation}
Then, under Assumptions \ref{ass:assnot2}.\ref{ass:Zij2} and \ref{ass:alpha0inf2}, there exists a constant $c>0$ depending on $\tau$, $||\alpha_0||_{\infty,\tau}$ and $\mathbb{E}[\e^{\bs{\beta_0^TZ}}]$ such that 
\[
\mathbb{P}((\widetilde\Omega^k_H)^c)\leq cn^{-k}.
\]
\end{lemma}
The proof of this lemma follows.
%\textcolor{red}{We refer to Subsection \ref{proof:lemHuang2} for the proof of this lemma. From Lemma \ref{lem:OmegaHann}, we obtain immediately the following corollary of Theorem \ref{th:Huangann2}.}
From this lemma, we can rewrite Theorem \ref{th:Huangann2} as:
\begin{corollary}
\label{cor:Huangann}
Let $\nu=B(\xi+1)s\Gamma_{n,k}/\{2\kappa^2(\xi,\mathcal{O})\}$, $k>0$ and $c>0$. Suppose Assumptions  \ref{ass:assnot2}.\ref{ass:Zij2} and \ref{ass:alpha0inf2} hold and $\nu\leq 1/\e$. Then, with probability larger than $1-cn^{-k}$
\[
|\bs{\hat\beta}-\bs{\beta_0}|_1\leq \dfrac{\e^{\eta}(\xi+1)s}{2\kappa^2(\xi,\mathcal{O})}\Gamma_{n,k} \quad \mbox{with} \quad \Gamma_{n,k}=C_0B\dfrac{\xi+1}{\xi-1}\sqrt{2\dfrac{\log(pn^k)}{n}},
\] 
where $\eta\leq 1$ is the smaller solution of $\eta\e^{-\eta}=\nu$ and $C_0>\sqrt{\tau||\alpha_0||_{\infty,\tau}\E[\e^{\bs{\beta_0^TZ}}]}$.
\end{corollary}

From Corollary \ref{cor:Huangann} and Assumption \ref{ass:assnot2}.\ref{ass:Zij2}, we deduce a prediction inequality given by the following proposition.
\begin{proposition}
\label{prop:predbetaann}
Let $k>0$ and $c>0$. Under Assumptions \ref{ass:assnot2}.\ref{ass:Zij2} and \ref{ass:assnot2}.\ref{ass:alpha0inf2}, with probability larger than $1-cn^{-k}$, we have
\begin{equation}
\label{predictionbetaann}
|\bs{\hat\beta}-\bs{\beta_0}|_1\leq C(s)\sqrt{\dfrac{\log (pn^k)}{n}},
\end{equation}
where $C(s)>0$ is a constant depending on the sparsity index $s$.
\newline
\end{proposition}

\begin{remark}
From Proposition \ref{prop:predbetaann} and Definition (\ref{set:OmegaH2}) of $\Omega^k_{H}$, we deduce that $\widetilde\Omega^k_H\subset\Omega^k_H$.
\end{remark}

%From Corollary \ref{cor:Huangann} and Assumption \ref{ass:Z}, we deduce a prediction inequality given by the following proposition.
%\begin{proposition}
%\label{prop:predbetaann}
%Let $k>0$ and $c>0$. Under Assumptions \ref{ass:Z} and \ref{ass:baseline3}.\ref{ass:alpha0inf3}, with probability larger than $1-cn^{-k}$, we have
%\begin{equation}
%\label{predictionbetaann}
%|\bs{\hat\beta^TZ}-\bs{\beta_0^TZ}|_1\leq C(s)\sqrt{\dfrac{\log (pn^k)}{n}},
%\end{equation}
%where $C(s)>0$ is a constant depending on the sparsity index $s$.
%\newline
%\end{proposition}
%
%\begin{remark}
%From Proposition \ref{prop:predbetaann} and Definition (\ref{set:OmegaH2}) of $\Omega^k_{H}$, we deduce that $\widetilde\Omega^k_H\subset\Omega^k_H$.
%\end{remark}
%

\textbf{Proof of Lemma \ref{lem:OmegaHann}}
\label{proof:lemHuang2}
To prove Lemma \ref{lem:OmegaHann}, we start from Lemma 3.3. p.10 in the paper of  \citet{Huang2013}, that we enounce below.

\begin{lemma}[Lemma 3.3 from \citet{Huang2013}]
\label{lem:probaHuang}
Suppose that Assumption \ref{ass:assnot2}.\ref{ass:Zij2} is verified. Let $\bs{\dot l^*_n(\bs{\beta})}$ be the gradient of the $ l^*_n(\bs{\beta})$ defined by (\ref{eq:PartialCoxLikelihood2}). Then, for all $C_0>0$,
\begin{equation}
\label{eq:lemHuang}
\mathbb{P}\left(|\bs{\dot l^*_n(\bs{\beta_0})}|_{\infty}>C_0Bx, \sumin\inttau Y_i(t)\diff N_i(t)\leq C^2_0n\right)\leq 2p\e^{-nx^2/2}.
\end{equation}
In particular, if $\max_{i\leq n}N_i(\tau)\leq 1$, then $\mathbb{P}(|\bs{\dot l^*_n(\bs{\beta_0})}|_{\infty}>Bx)\leq 2p\e^{-nx^2/2}$.
\end{lemma}

Before proving the lemma that is in interest, we recall the Bernstein Inequality for martingales (see \citet{SVG}).

\begin{lemma}[Lemma 2.1 from \citet{SVG}]
\label{lem:svg}
Let $\{M_t\}_{t\geq 0}$ be a locally square integrable martingale w.r.t. the filtration $\{\mathcal{F}_t\}_{t\geq 0}$. Denote the predictable variation of $\{M_t\}$ by $V_t=\langle M,M \rangle_t$, $t\geq 0$, and its jumps by $\Delta M_t=M_t-M_{t^-}$.
Suppose that $|\Delta M(t)|\leq K$ for all $t>0$ and some $0\leq K<\infty$. Then for each $a>0$, $b>0$,
\[
\mathbb{P}(M_t\geq a\mbox{ and } V_t\leq b^2 \mbox{ for some t })\leq \exp\left[-\dfrac{a^2}{2(aK+b^2)}\right].
\]
\end{lemma}

From Lemma \ref{lem:probaHuang}, to prove Lemma \ref{lem:OmegaHann}, it remains to control
\[
\mathbb{P}\left(\sumin\inttau Y_i(t)\diff N_i(t)> C^2_0n\right),
\]
%|\bs{\dot l^*_n(\bs{\beta_0})}|_{\infty}>C_0Bx, 
%to remove the event $\{\sum_{i=1}^{n}\int_{0}^{\tau} Y_i(t)\diff N_i(t)> C^2_0n\}$ from the probability (\ref{eq:lemHuang}) and prove Lemma \ref{lem:OmegaHann}.
Using the Doob-Meyer decomposition  and since, 
\[
\sum_{i=1}^{n}\int_{0}^{\tau}Y_i(t)\alpha_0(t)\expbeta Y_i(t)\diff t\leq n\tau||\alpha_0||_{\infty,\tau}\e^{B|\bs{\beta_0}|_1},
\]
we obtain for $C_0>\sqrt{\tau||\alpha_0||_{\infty,\tau}\E[\e^{\bs{\beta_0^TZ}}]}$,
\begin{align*}
\mathbb{P}\Bigg(\sumin\inttau Y_i(t)\diff N_i(t)>C_0^2n\Bigg)
\leq \mathbb{P}\Bigg(\sumin\inttau Y_i(t)\diff M_i(t)>C^2_0n-n\tau||\alpha_0||_{\infty,\tau}\e^{B|\bs{\beta_0}|_1}\Bigg).
\end{align*}
Then, we apply Lemma \ref{lem:svg} to the martingale $\sum_{i=1}^{n}\int_{0}^{\tau} Y_i(t)\diff M_i(t)$, with $K=1$ and 
\begin{align*}
V_t=\mathbb{E}\Big[\sumin\inttau Y^2_i(t)\alpha_0(t)\expbeta Y_i(t)\diff t\Big]\leq ||\alpha_0||_{\infty,\tau}\tau\mathbb{E}[\e^{\bs{\beta_0^TZ}}]n.
\end{align*}
We obtain
\begin{align*}
\mathbb{P}\Bigg(\sumin\inttau Y_i(t)\diff M_i(t)&>C^2_0n-n\tau ||\alpha_0||_{\infty,\tau}\mathbb{E}[\e^{\bs{\beta_0^TZ}}]\Bigg)\\
&\leq \exp\Bigg(-\dfrac{n(C^2_0-\tau||\alpha_0||_{\infty,\tau}\mathbb{E}[\e^{\bs{\beta_0^TZ}}])^2}{2C^2_0}\Bigg).
\end{align*}
Finally, we get
\begin{align*}
\mathbb{P}\Big(|\bs{\dot l^*_n(\bs{\beta_0})}|_{\infty}>C_0Bx\Big)\leq 2p\e^{-nx^2/2}+\exp\Bigg(-\dfrac{n}{2C^2_0}(C^2_0-\tau||\alpha_0||_{\infty,\tau}\mathbb{E}[\e^{\bs{\beta_0^TZ}}])\Bigg).
\end{align*}

Taking $x=\sqrt{2\log(n^kp)/n}$, there exists a constant $c>0$ depending on $\tau$, $||\alpha_0||_{\infty,\tau}$ and $\mathbb{E}[\e^{\bs{\beta_0^TZ}}]$ such that 
\[
\mathbb{P}((\widetilde\Omega^k_H)^c)\leq cn^{-k},
\]
which leads to the expected result of Lemma \ref{lem:OmegaHann}.
\qed

%
%
%
%
%\newpage
\bibliography{biblio3}
\bibliographystyle{plainnatfrench}
\end{document}